\documentclass[journal ]{new-aiaa}
\usepackage[utf8]{inputenc}
\usepackage{textcomp}

\usepackage{graphicx}
\usepackage{amsmath}
\usepackage[version=4]{mhchem}
\usepackage{siunitx}
\usepackage{longtable,tabularx}
\usepackage{bbm}
\usepackage{algorithm}
\usepackage{algpseudocode}
\usepackage{multirow}

\title{Spare Strategy Analysis and Design\\ for Mega Satellite Constellations Using Markov Chain}

\author{Seungyeop Han \footnote{Ph.D. Candidate, Daniel Guggenheim School of Aerospace Engineering}}
\author{Zachary Grieser \footnote{Ph.D. Student, Daniel Guggenheim School of Aerospace Engineering}}
\affil{Georgia Institute of Technology, Atlanta, Georgia, 30332}
\author{Shoji Yoshikawa\footnote{Chief Expert, Advanced Technology R\&D Center, Mitsubishi Electric Corporation}}
\author{Takumi Noro\footnote{Researcher, Advanced Technology R\&D Center, Mitsubishi Electric Corporation}}
\author{Takumi Suda\footnote{Head Researcher, Advanced Technology R\&D Center, Mitsubishi Electric Corporation}}
\affil{Mitsubishi Electric Corporation, Amagasaki 661-8861, Japan}
\author{Koki Ho \footnote{Dutton-Ducoffe Professor, Associate Professor, Daniel Guggenheim School of Aerospace Engineering, AIAA Senior Member, kokiho@gatech.edu (Corresponding Author)}}
\affil{Georgia Institute of Technology, Atlanta, Georgia, 30332}

\begin{document}
\footnotetext{This paper is a substantially revised version of a portion of Paper AAS 24-165 presented at the AAS/AIAA Astrodynamics Specialist Conference, Broomfield, CO, August 11-15, 2024.}
\maketitle

\begin{abstract}
This paper presents a Markov-chain-based method for the early-phase analysis and design of spare-management architectures for large-scale satellite constellations. To assess the long-run viability of such concepts of operations, satellite failure and replenishment processes are modeled as Markov chains and analyzed through their stationary solution. We reinvestigate an indirect spare strategy, modeled as a multi-echelon periodic-review reorder-point/order-quantity policy, in which spares are first delivered to parking orbits and then transferred to constellation planes.
The stock levels in constellation and parking orbits are each modeled as independent Markov chains, and a fixed-point iteration yields a consistent joint stationary solution that describes the strategy’s average behavior. This approach accurately captures the stochastic interplay within a multi-echelon model driven by orbital mechanics, avoiding the aggregation assumptions of prior works and remaining valid across a wider operating domain. Building on this fast, accurate analysis, we formulate an optimization problem and solve it via a genetic algorithm. Finally, we demonstrate the practical value of both the analysis method and the optimization framework in a real-world mega-constellation case study.
\end{abstract}

\section*{Nomenclature}
{\renewcommand\arraystretch{1.0}
\noindent\begin{longtable*}{@{}l @{\quad=\quad} l@{}}
$\tau_\text{mc}$ & Time step of the discrete-time Markov process, in days \\
$\lambda_\text{sat}$ & Failure rate of a satellite, in failures per unit time $\tau_\text{mc}$\\
$\mu_\text{lv}$ & Mean interval between launches, in days \\
$\tau_\text{lv}$ & Constant launch order processing time, in days \\
$\tau_\text{c}/\tau_\text{p}$ & RAAN alignment period from the perspective of in-plane/parking orbits \\
$q_{\text{c}}/q_{\text{p}}$ & Replenishment quantity for in-plane/parking spares \\
$r_{\text{c}}/r_{\text{p}}$ & Reorder point for in-plane/parking spares \\
$ N_{\text{sat}_\text{c}}/N_{\text{sat}_\text{p}}$ & 
Maximum in-plane/parking state level including operational and spare satellites \\
$\bar{N}_\text{sat}$ & Nominal number of operational satellites per in-plane orbit \\
$ N_{\text{orbit}_{\text{c}}} / N_{\text{orbit}_{\text{p}}}$  & Number of in-plane (constellation)/parking orbital planes \\
$ P_{f_{\text{c}}}/P_{f_{\text{p}}}$ & Failure transition matrix for in-plane/parking states per review period \\
$ P_{q_{\text{c}}}/ P_{q_{\text{p}}}$ & $q$-unit replenishment transition matrix for in-plane/parking states\\
$\pi^{q_{\text{c}}}/\pi^{q_{\text{p}}}$ & Expected in-plane/parking state distribution immediately after $q$-unit replenishment \\
$\pi^{r_{\text{c}}}/\pi^{r_{\text{p}}}$ & Expected in-plane/parking state distribution at the $r$-reorder point \\
$\pi^{\text{io}_{\text{c}}}/\pi^{\text{io}_{\text{p}}}$ & Expected in-plane/parking state distribution during the inter-order (IO) period \\
$\pi^{\text{lt}_{\text{c}}}/\pi^{\text{lt}_{\text{p}}}$ & Expected in-plane/parking state distribution during the lead-time (LT) period \\
$\pi^{\text{rc}_{\text{c}}}/\pi^{\text{rc}_{\text{p}}}$ & Expected in-plane/parking state distribution over the full replenishment cycle (RC) \\
\end{longtable*}}

\section{Introduction}
\lettrine{T}{he} satellite constellation is a coordinated group of satellites that together enable continuous global coverage. Constellations with thousands of satellites are becoming central to space-based communication infrastructure. Starlink has deployed over 9{,}000 satellites as of 2025 and plans to expand to more than 42{,}000~\cite{starlink1,starlink2}. Amazon’s Project Kuiper began launching in 2025 and targets 3{,}236 satellites~\cite{kuiper2025}. OneWeb, now under Eutelsat, has completed its initial deployment of 648 satellites and is planning expansion~\cite{oneweb}. China’s Guowang network is expected to include about 13{,}000 satellites and has already begun initial deployment~\cite{guowang2025}. These systems rely on advanced constellation architecture, frequency coordination, and inter-satellite communication, with a comparative analysis of their designs provided by Pachler et al.~\cite{Pachler2021Comp}.

Over time, satellite failures reduce operational capacity and degrade constellation performance. To maintain service, an effective population management strategy is required. Two representative approaches are spare deployment and on-orbit servicing. The former inserts a spare satellite to replace a failed one~\cite{cornara1999satellite}, while the latter uses a servicer to repair the failed satellite~\cite{luu2022orbit}.  Spare deployment is generally more suitable for mega-constellations with small, low-cost satellites, whereas servicing is more appropriate for smaller systems with high-value assets. In the context of low Earth orbit (LEO) mega-constellations, launching replacement satellites is widely regarded as more cost-effective than performing on-orbit repairs. This is driven by several factors: declining launch costs~\cite{potter2023mooreslaw}, economies of scale from mass production, and the high expense of designing satellites to be serviceable. Additionally, small satellites can have reliability challenges related to limited testing and a lack of redundancy~\cite{BOUWMEESTER2022108288}, further reinforcing the need for rapid spare deployment. These trends motivate the study of optimal spare-management strategies, which is the focus of this paper.

Early work on spare deployment modeled satellite replenishment using an $(s,S)$ inventory framework, where $s$ is the reorder point and $S$ is the order-up-to level, with exponential lead times, establishing a probabilistic basis for maintaining constellation performance~\cite{1966dishon}. Subsequent research analyzed various spare deployment strategies—including in-orbit, parking orbit, and ground-based options—across different phases of the constellation life cycle~\cite{cornara1999satellite}. Discrete-time Markov decision processes were later applied to determine cost-optimal satellite replacement policies over finite planning horizons~\cite{sumter2003optimal}. A more scalable formulation for mega-constellations was introduced through a multi-echelon reorder-point/order-quantity policy, commonly referred to as the $(s,Q)$ inventory model, or equivalently the $(r,q)$ model, that considers parking orbits as intermediate warehouses between the ground and constellation planes~\cite{jakob2019optimal}. Here, $s$ (or $r$) is the reorder point and $Q$ (or $q$) is the order quantity.

Beyond inventory modeling, other studies have explored resilience and value-based assessment. A value model was developed to examine trade-offs between satellite reliability and system attributes such as cost or mass~\cite{collopy2003assigning}. A quantitative resilience framework was proposed, incorporating robustness and recovery metrics to support mission-driven constellation optimization~\cite{cuhran2017resilience}. Resilience under solar weather radiation was also analyzed using a time-based simulation approach that evaluated reconstitution strategies and highlighted the influence of manufacturing and launch timelines on system performance~\cite{novak2023analysis}. However, while prior work has studied satellite spare-replenishment strategies and separate studies have examined constellation value or resilience, the literature still lacks an accurate and computationally efficient framework for early-phase assessment of indirect multi-echelon spare architectures with coupled parking-orbit availability, in-plane demand, and replenishment timing.

Building on previous work~\cite{han2024analysis}, this paper develops a Markov-based framework for indirect multi-echelon spare architectures in large-scale satellite constellations. The approach employs Markov processes to represent state transitions and analyze long-run behavior through stationary state distributions. The framework incorporates key physical realities often simplified in other models, such as the non-negativity of physical stock and state-dependent failure rates. The central contribution is the high-fidelity modeling of the multi-echelon replenishment dynamics. Unlike prior works that use aggregation assumptions, our approach models the stochastic interplay between the probability distribution of in-plane demand and that of parking orbit spare availability, a dynamic physically constrained by the periodic lead times from orbital mechanics. By computing the full state probability distribution, a step often intractable in prior analytical models, the framework enables a rigorous analysis of operational costs and constellation availability. We demonstrate the framework on a real-world mega-constellation case study and illustrate its use in early-phase architecture assessment and spare-strategy design.

The remainder of the paper is organized as follows. Section~\ref{sec2} introduces the modeling preliminaries. Section~\ref{sec3} presents the proposed analytical method for evaluating the indirect resupply strategy. Section~\ref{sec4} applies this method to assess system performance. Section~\ref{sec5} validates the model through comparison with Monte Carlo simulations, and Section~\ref{sec6} illustrates its application in a design optimization context. Finally, Section~\ref{sec7} concludes the paper and suggests directions for future research.

\section{Preliminaries} \label{sec2}
This section introduces the spare-management architecture and the physical and probabilistic ingredients needed to formulate the stationary Markov-chain model for the indirect strategy in Sec.~\ref{sec3}.
\subsection{Spare Management Policy}
\subsubsection{Indirect Resupply Strategies}

The indirect strategy uses a large launch vehicle (LV) to deliver spares to a parking orbit, from which they are transferred to the constellation orbit via drift in  the right ascension of the ascending node (RAAN). In this architecture, a shared parking-orbit layer, consisting of $N_{\text{orbit}_\text{p}}$ parking orbits, serves all $N_{\text{orbit}_\text{c}}$ constellation planes. Although this method benefits from batch discounts and lower launch cost, the RAAN drift introduces a longer replenishment delay. By contrast, the direct strategy employs a small LV to replenish the in-plane orbit immediately. A detailed explanation of the direct strategy can be found in past work~\cite{han2024analysis}.

\begin{figure}[!h]
    \centering
    \includegraphics[width=.45\textwidth]{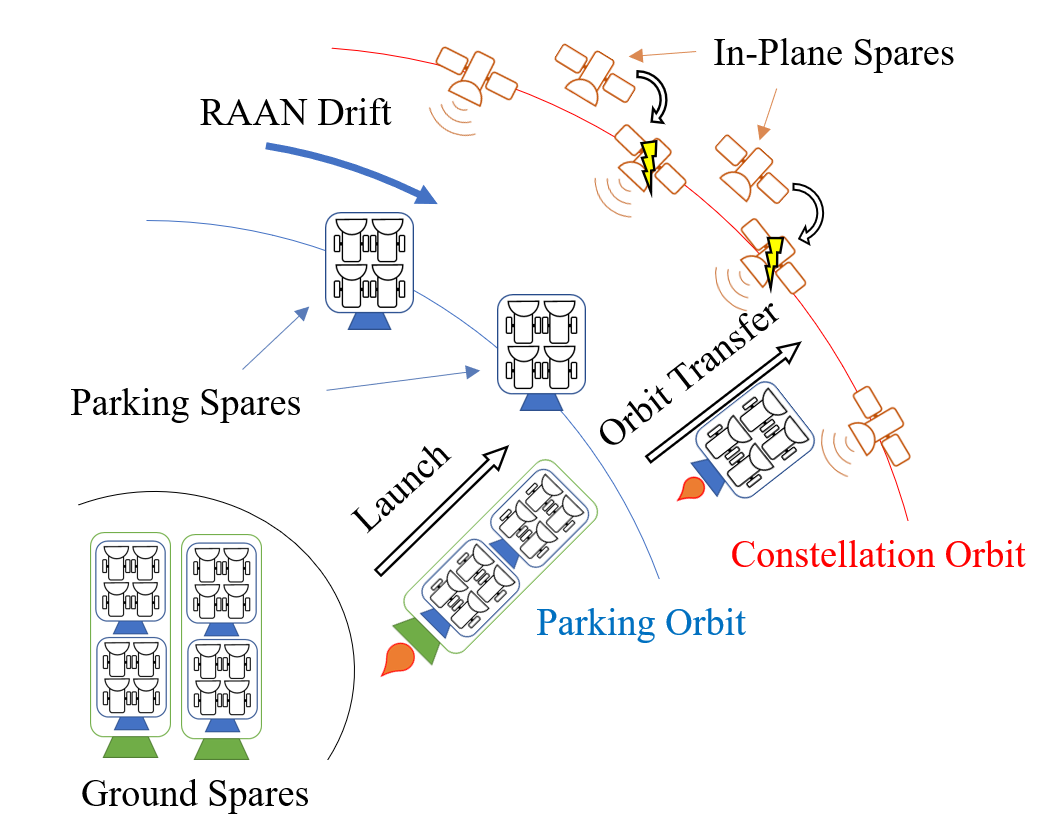}
    \caption{Illustration of Indirect Spare Strategy}
    \label{fig_strategy}
\end{figure}

Figure~\ref{fig_strategy} illustrates the indirect strategy. If a failure occurs, an in-plane spare immediately replaces the failed satellite, and whenever the number of in-plane spares falls below a threshold, a batch of spares is transferred from one of the parking orbits. At the same time, if the number of parking spares falls below a second threshold, a ground resupply order is placed and the replacement arrives after the LV’s lead time. To enable these transfers, spare satellites are grouped into transfer buses and launched in batches aboard the large LV, allowing spares to be delivered to multiple constellation planes as needed.

\subsubsection{Inventory Management Policy} \label{sec:inventory_model} 
\begin{figure}[!ht]
    \centering
    \includegraphics[width=.45\textwidth]{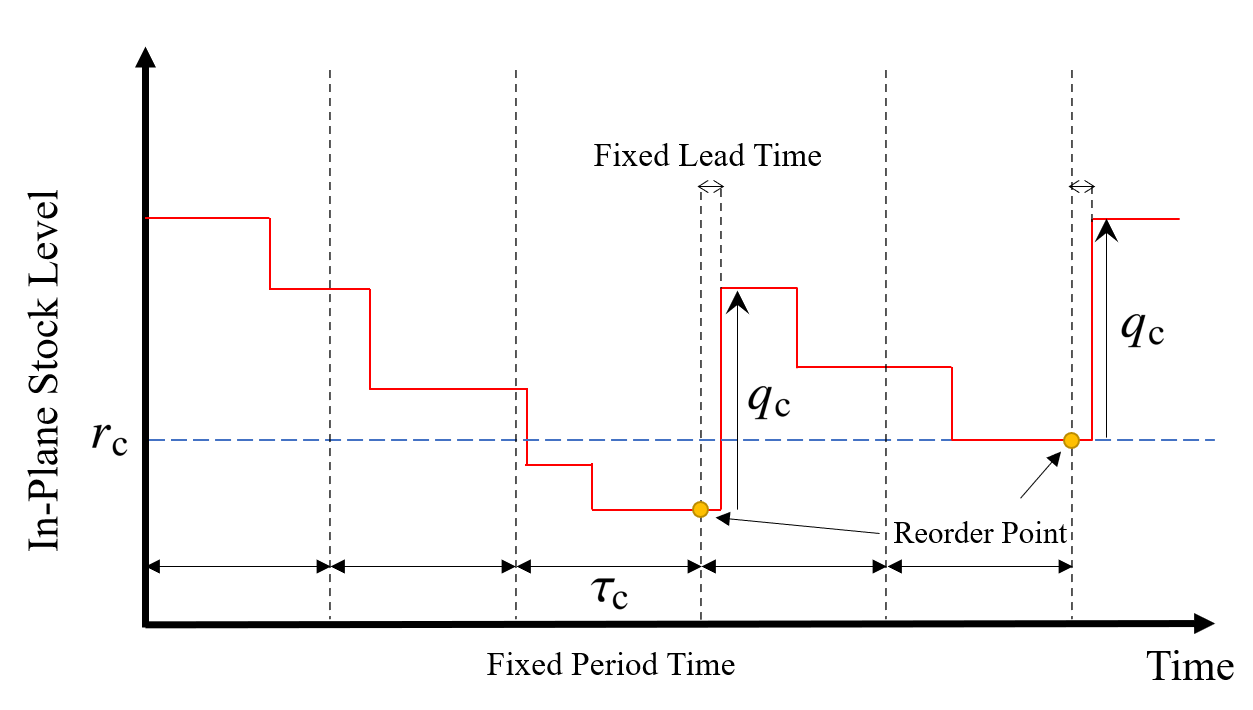}
    \includegraphics[width=.45\textwidth]{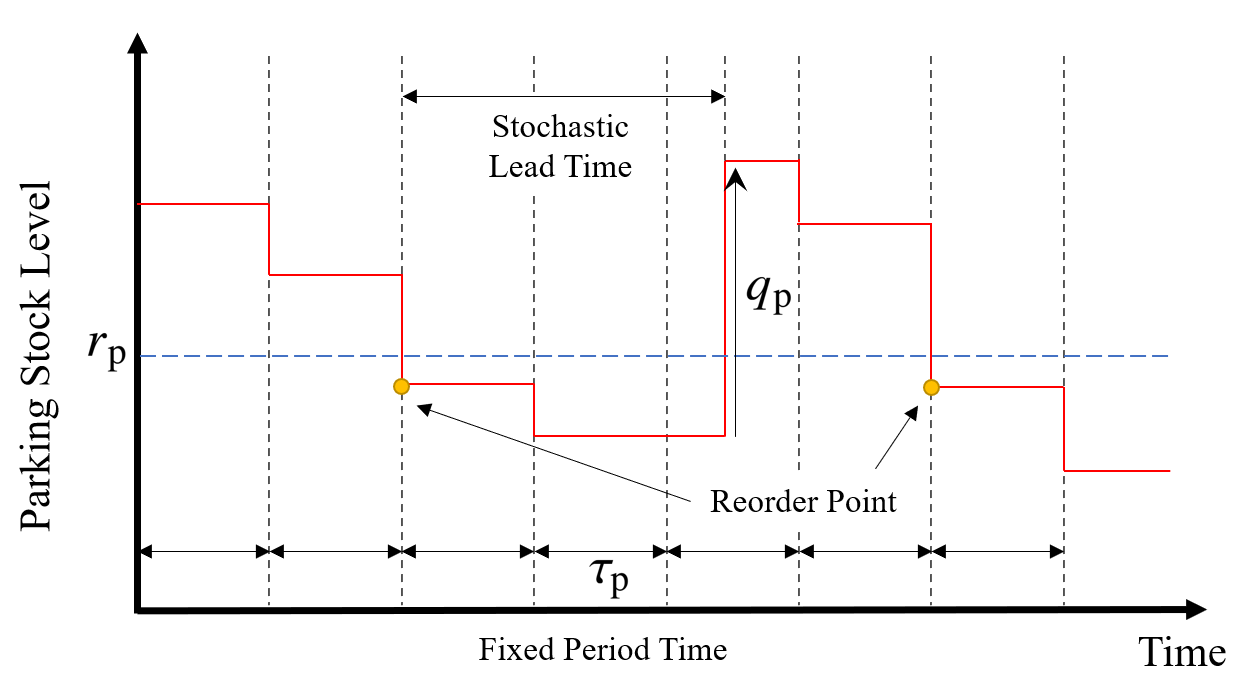}
    \caption{Stock level profile of (a) constellation and (b) parking orbits under $(r, q, \tau)$ policy}
    \label{fig:inventory_model}
\end{figure}

For the indirect method, both constellation and parking orbits follow the same $(r,q,\tau)$ policy structure, where $r$ is the reorder point, $q$ is the order quantity, and $\tau$ is the review period, but with different review-period and lead-time parameters. In the in-plane orbit, the stock is reviewed every $\tau_\text{c}$ and, if $X_\text{c}\le r_\text{c}$, an order of size $q_\text{c}$ is placed and arrives after the fixed orbital-transfer lead time; in the parking orbit, the stock is reviewed every $\tau_\text{p}$ and, if $X_\text{p}\le r_\text{p}$, an order of size $q_\text{p}$ is placed and arrives after the stochastic LV lead time. Figure~\ref{fig:inventory_model}(a) and (b) illustrate these two cases, respectively.

\subsection{Orbital Mechanics}
\subsubsection{Constellation Model}
This research focuses on large-scale constellations in low Earth orbit (LEO), particularly the well-known Walker Delta pattern constellation~\cite{walker1984satellite}. This constellation comprises multiple circular constellation orbits with identical inclination angles. Additionally, 
their RAAN values $\Omega$ are evenly distributed. In this configuration, there are a total of $N_{\text{orbit}_\text{c}}$ constellation orbital planes (also referred to as in-plane orbits), with nominal $\bar{N}_{\text{sat}}$ satellites allocated to each constellation plane.

\subsubsection{Parking Orbits}
The parking orbits are temporary orbits where the batch of spares is first inserted for the indirect resupply strategy. There are a total of $N_{\text{orbit}_\text{p}}$ parking orbits, and they are evenly distributed in RAAN. The parking orbits are assumed to have the same inclination angle as the constellation orbits but have different semi-major axes.

\subsubsection{Orbital Maneuver}
This research assumes that all orbit transfers are conducted using coplanar Hohmann transfers, and out-of-plane maneuvers are excluded due to their inefficiency in terms of delta-v. To facilitate these transfers, spare satellites are grouped with a transfer bus into a single batch. 
For a propulsive maneuver, the mass of fuel $m_\text{fuel}$ required for a given $\Delta v$ is computed as:
\begin{equation} \label{eq:fuel_mass}
    m_\text{fuel} = m_\text{dry} (e^{\Delta v / v_\text{ex}}- 1 )
\end{equation}
where $m_\text{dry}$ is the dry mass of the spacecraft (e.g., set of spares and transfer bus), $v_\text{ex}$ is the effective exhaust velocity and $\Delta v$ is the velocity increment for the transfer. For the Hohmann transfer considered here, the required $\Delta v$ is calculated as follows:
\begin{equation} \label{eq:delta_v}
    \Delta v = \sqrt{\frac{\mu_\oplus}{a_\text{p} }}
    \left( \sqrt{\frac{2a_\text{c} }{a_\text{p} +a_\text{c} }} - 1\right)
    + \sqrt{\frac{\mu_\oplus}{a_\text{c} }}
    \left( 1-\sqrt{\frac{2a_\text{p} }{a_\text{p} +a_\text{c} }} \right)
\end{equation}
where $\mu_\oplus$ is the Earth's gravitational constant, and $a_\text{p} $ and $a_\text{c} $ ($a_\text{c} >a_\text{p} $) are the radii of initial (parking) and final (constellation) orbits, respectively.

The time of flight for a Hohmann transfer in LEO is typically on the order of a few hours. As this duration is negligible compared to the time horizon of the management policy, the maneuver is assumed to be instantaneous for modeling simplicity.

\subsubsection{RAAN Drift}
The LEO experiences non-negligible perturbations due to the oblateness of the Earth, resulting in the secular drift of the orbital plane. With an assumption of circular orbits, the RAAN drift rate of constellation and parking orbits are computed as:
\begin{equation} \label{eq:rann_drift}
    \dot \Omega_\text{c} = -\frac{3}{2} J_2 \frac{R_\oplus^2 }{a_\text{c}^2} \sqrt{\frac{\mu_\oplus}{a_\text{c}^3}} \cos i,\quad \dot \Omega_\text{p} = -\frac{3}{2} J_2 \frac{R_\oplus^2 }{a_\text{p}^2} \sqrt{\frac{\mu_\oplus}{a_\text{p}^3}}  \cos i
\end{equation}
where $a_\text{c}$ and $a_\text{p}$ are the semi-major axis, $i$ is the inclination of the constellation and parking orbit, $R_\oplus$ is the Earth's radius, and $J_2$ is the second zonal harmonic coefficient of the Earth\cite{prussing1993orbital}.

\subsection{Markov Chain Model}
We model the spare‐satellite count as a discrete‐time Markov chain.  Let \(X_k\in\{0,1,\dots,N_{\rm sat}\}\) be the number of satellites (including spares) at step \(k\), and write the distribution

\begin{equation} 
    \pi_k = \begin{bmatrix}
    \mathbb{P}(X_k = N_\text{sat}) & \mathbb{P}(X_k = N_\text{sat}-1) & \cdots & \mathbb{P}(X_k = 0)
\end{bmatrix}^\top
\end{equation}
Here, the entries of $\pi_k$ are ordered from state $N_\text{sat}$ down to state $0$.
Assuming a time‐homogeneous transition matrix $P\in\mathbb R^{(N_{\rm sat}+1)\times(N_{\rm sat}+1)}$ with entries $P_{ij}=\mathbb{P}(X_{k+1}=i\mid X_k=j)$, the chain evolves by $\pi_{k+1}=P\,\pi_k$. Under the usual ergodicity conditions (e.g.\ $P$ irreducible and aperiodic), the Markov Chain has a unique $\pi$ satisfying
\begin{equation} 
    \pi = P \pi
\end{equation}
and it gives the long‐run fraction of time the chain spends in each state\cite{hillier2005introduction}. In our spare‐management model, failures (state decreases) and replenishment (state increases) guarantee these conditions.

\subsection{Probabilistic Model}
\subsubsection{Satellite Failure Probability Distribution}
Satellite failures are modeled as a Poisson process (i.e. exponential interarrival times), and spare satellites are assumed not to fail. 
Although more accurate lifetime laws (e.g., Weibull or bathtub models) exist \cite{castet2009satellite}, the present paper focuses on infinite-horizon stationary analysis for early-phase assessment rather than transient age-dependent behavior. Under the assumed spare policy, failed satellites are immediately replaced by in-orbit spares and the spare stock is then replenished through inventory control. Accordingly, the long-run behavior is represented by an effective constant failure-rate model, making the Poisson approximation reasonable.

Let $\tau_\text{mc}$ denote the time step of the Markov process, and let $\lambda_\text{sat}$ be the failure rate of an operational satellite per $\tau_\text{mc}$. Then, the probability of observing $k$ failures from $n$ satellites (including spares) during $\tau_\text{mc}$ is given by:
\begin{equation}\label{mpois_fail}
    \nu_{k,n} = 
    \mathbb{P}(F=k| X=n) = \begin{cases}
        0 & \mbox{if } k > \bar{N}_\text{sat}\\ 
        \frac{(n\lambda_\text{sat})^k}{k!}e^{-n\lambda_\text{sat}} & \mbox{if } n \leq \bar{N}_\text{sat} \mbox{ and } k \leq \bar{N}_\text{sat} \\ 
        \frac{(\bar{N}_\text{sat}\lambda_\text{sat})^k}{k!}e^{-\bar{N}_\text{sat}\lambda_\text{sat}} & \mbox{if } n > \bar{N}_\text{sat} \mbox{ and } k \leq \bar{N}_\text{sat} \\
    \end{cases}
\end{equation}
where $F$ is the number of failures, $\bar{N}_\text{sat}$ is the nominal number of operational (non-spare) satellites. This formulation assumes immediate failure replacement and that spare satellites do not fail (i.e., $k>\bar{N}_\text{sat}$). When $n \leq \bar{N}_\text{sat}$, all satellites are operational, yielding a total failure rate of $n\lambda_\text{sat}$. When $n > \bar{N}_\text{sat}$, the number of operational satellites is fixed at $ \bar{N}_\text{sat}$, making a total failure rate of $\bar{N}_\text{sat} \lambda_\text{sat}$, and excess satellites are treated as spares.

Additionally, let $N_{\text{sat}_\text{c}}$ be the maximum number of satellites including spares in constellation orbits. Then the state transition matrix due to failure can be defined as:
\begin{equation} \label{eq:fail_matrix_sim}
    P_{f} = \begin{bmatrix}
    \nu_{0,N_{\text{sat}_\text{c}}} 
      & 0                     & \cdots         & 0         \\
    \nu_{1,N_{\text{sat}_\text{c}}} 
      & \nu_{0,N_{\text{sat}_\text{c}}-1}   & \cdots         & 0         \\
    \vdots 
      & \vdots                & \ddots & \vdots       \\
    1-\sum_{k=0}^{N_{\text{sat}_\text{c}}-1} \nu_{k,N_{\text{sat}_\text{c}}}
      & 1-\sum_{k=0}^{N_{\text{sat}_\text{c}}-2} \nu_{k,N_{\text{sat}_\text{c}}-1} 
                             & \cdots 
                                       & \nu_{0,0}  
\end{bmatrix}
\end{equation}
where $P_f \in \mathbb{R}^{( N_{\text{sat}_\text{c}}+1 )\times ( N_{\text{sat}_\text{c}}+1)}$. 
Here, for rows above the last row, element $(a,b)$ is $\nu_{a-b,\,N_{\text{sat}_\text{c}}-b+1}$, i.e., the probability of $a-b$ failures from state $X=N_{\text{sat}_\text{c}}-b+1$. The last row collects all probability mass associated with transitions to the lowest state, ensuring that each column of $P_f$ sums to one. Since the states are ordered from high to low stock level, the lower-triangular structure of $P_f$ means that probability mass can only remain in the same state or move to a lower stock state. In summary, if $\pi$ is the in-plane state distribution, then $P_f\pi$ gives the distribution after a one-step failure.

\subsubsection{LV Lead‐Time Probability Distribution}
The ground‐resupply lead time is modeled as a shifted exponential distribution \cite{jakob2019optimal}: $T \sim \text{Exp}(\mu_\text{lv}) + \tau_\text{lv}$, and its probability density function is
\begin{equation}
    f(T=t;\mu_\text{lv},\tau_\text{lv}) = 
    \begin{cases}
        \frac{1}{\mu_\text{lv}} e^{-{(t-\tau_\text{lv})}/{\mu_\text{lv}}} & \quad t \geq \tau_\text{lv}\\
        0 & \quad t < \tau_\text{lv}\\
    \end{cases}
\end{equation}
where $\mu_\text{lv}$ is the mean of the exponential component and $\tau_\text{lv}$ is the  fixed LV‐processing delay. To simplify our discrete‐time modeling, we choose $\tau_\text{mc}$ such that $\tau_\text{lv}$ is an integer multiple of $\tau_\text{mc}$, i.e., $\tau_\text{lv} = k_\text{lv} \tau_\text{mc}$. Then the probability of having a lead time between $k$ and $k+1$ time steps of $\tau_\text{mc}$ is computed as
\begin{equation} \label{eq:rho_lv}
\begin{aligned}
    \rho_{k+1} &= \mathbb{P}(k\tau_\text{mc} \leq T < (k+1) \tau_\text{mc}) \\
    &= 
    \begin{cases}
        e^{(-k\tau_\text{mc} + \tau_\text{lv})/{\mu_\text{lv}}}\left( 1 - e^{-\tau_\text{mc}/{\mu_\text{lv}}}\right), & \text{if } k\tau_\text{mc} \geq \tau_\text{lv} \\
        0, & \text{otherwise } \\
    \end{cases}
\end{aligned}
\end{equation}
Note that each orbit is assumed to place at most one LV order at a time.

\section{Modeling and Analysis of Spare Management Policy} \label{sec3}
The indirect resupply strategy is analyzed by solving the constellation‐orbit and parking‐orbit Markov chains independently, then enforcing consistency between them. Specifically, the constellation stock level depends on parking‐orbit spare availability, while the parking stock level is driven by demand from the constellation orbits. In other words, the two chains are coupled through in‐plane spares demand and parking spares availability. In the next subsection, we present each analysis step by step, and at the end of this section we show how to combine them.

Lastly, the subscripts $(\cdot)_\text{c}$ and $(\cdot)_\text{p}$ denote the constellation and parking orbits, respectively. The maximum number of satellites in constellation orbit is $ N_{\text{sat}_\text{c}} = q_\text{c} + r_\text{c}$ in units of satellites, and the maximum number of parking spares is $ N_{\text{sat}_\text{p}} = q_\text{p} + r_\text{p}$ in units of batches. That is, $\pi^{(\cdot)_\text{c}} \in \mathbb{R}^{{N}_{\text{sat}_\text{c}}+1}$ and $\pi^{(\cdot)_\text{p}} \in \mathbb{R}^{ N_{\text{sat}_\text{p}}+1}$. 

\subsection{Repeating Structure of In-plane and Parking Orbits}
As previously explained, the indirect strategy uses the relative RAAN drift between constellation and parking orbits to transfer spare satellites. Because the orbital planes are symmetrically distributed, each parking orbit and each constellation orbit align periodically. From a constellation orbit’s perspective, the interval between successive alignments with parking orbits is
\begin{equation} \label{eq:t_plane_indirect}
    \tau_\text{c} = \frac{2\pi}{N_{\text{orbit}_\text{p}} \lvert \dot \Omega_\text{c}  - \dot \Omega_\text{p}\rvert }
\end{equation}
Conversely, from a parking orbit’s perspective, the interval between successive alignments with constellation orbits is
\begin{equation} \label{eq:t_park_indirect}
    \tau_\text{p} = \frac{2\pi}{N_{\text{orbit}_\text{c}} \lvert \dot \Omega_\text{c}  - \dot \Omega_\text{p}\rvert }
\end{equation}
With a proper selection of $\tau_\text{mc}$ or through rounding operations, $\tau_\text{c}$ and $\tau_\text{p}$ can be expressed as integer multiples of $\tau_\text{mc}$ as:
\begin{equation} 
    \tau_\text{p} = k_\text{p} \tau_\text{mc}, \quad 
    \tau_\text{c} = k_\text{c} \tau_\text{mc} \quad k_\text{p}, k_\text{c} \in \mathbb{N}
\end{equation}
Here, $\tau_c$ and $\tau_p$ are the actual review periods of the $(r,q,\tau)$ policy, determined by the relative RAAN drift, whereas $\tau_\text{mc}$ is only the numerical discretization used to represent those periods and the LV lead-time model.


\subsection{Constellation Orbit Analysis Method} \label{sec:in_plane_anal_indirect}
\subsubsection{Replenishment Transition Matrix} \label{sec:in_plane_replenshiment_indirct}
Assume each parking‐orbit stock \(X_\text{p}\) is independent and identically distributed, and let its state distribution just before a RAAN contact be known.  We then define the parking‐availability probability (see Section~\ref{sec:park_avail}) as
\begin{equation} \label{eq:parking_avail}
    \kappa_j = \mathbb{P}(X_\text{p} \geq j | E)
\end{equation}
where $E$ denotes the RAAN contact event. Here, $j$ denotes the realized demand level, i.e., $D_\text{c}=j$, where $D_\text{c}$ is the random variable for the demand in batches of size $q_\text{c}$ from the in-plane orbit:
\begin{equation} \label{eq:park_demand_func}
    D_\text{c} = \begin{cases}
        \Bigl\lceil\frac{ r_\text{c} + 1 - X_\text{c}}{q_\text{c}}\Bigr\rceil, & \text{if}\ X_\text{c} \leq r_\text{c}, \\
        0, & \text{if}\ X_\text{c} > r_\text{c},
    \end{cases}
\end{equation}
where $\lceil \cdot \rceil$ is the ceiling operator. For example, the demand is zero when the number of satellites exceeds $r_\text{c}$. On the other hand, if $X_\text{c}\leq r_\text{c}$, the demand corresponds to the number of batches required to raise the in-plane stock level above $r_\text{c}$. With these  definitions, \(\kappa_j\) is the probability that at least \(j\) batches are available in the parking orbit when RAANs align.

Then the in-plane replenishment transition matrix becomes:
\begin{equation} \label{eq:park_avail_mat}
    P_{q_\text{c}} = \begin{bmatrix}
        \kappa_0I_{q_\text{c}} & \kappa_1 I_{q_\text{c}} & \kappa_2 I_{q_\text{c}}  & \cdots \\
        0 & (\kappa_0 -\kappa_1)I_{q_\text{c}} & (\kappa_1 -\kappa_2) I_{q_\text{c}}  & \cdots \\
        0 & 0 & (\kappa_0 -\kappa_1) I_{q_\text{c}}  & \cdots \\
        \vdots & \vdots & \vdots  & \ddots
    \end{bmatrix},
\end{equation}
and $P_{q_\text{c}}\in \mathbb{R}^{({N}_{\text{sat}_\text{c}}+1)\times ({N}_{\text{sat}_\text{c}}+1)}$. Each \(q_c\times q_c\) block reflects the demand level and corresponding parking spare availability. For example, the first $q_\text{c}$ entries of $\pi$ (i.e., $r_\text{c}<X_\text{c}\leq N_{\text{sat}_\text{c}}$) require no replenishment so are multiplied by $\kappa_0 I_{\text{q}_\text{c}}$. The next $q_\text{c}$ entries (i.e., $r_\text{c}-q_\text{c}<X_\text{c}\leq r_\text{c}$) has a demand of one batch. If the parking spare is available with $\kappa_1$ probability (i.e., (1,2) block), the stock level increase by $q_\text{c}$ after replenishment. Otherwise, with probability $\kappa_0 - \kappa_1$ (i.e., (2,2) block), the state level remains unchanged. In summary, let $\pi$ be the distribution right before receiving parking spares, then $P_{q_\text{c}}\pi$ is the distribution after receiving parking spares considering the parking availability.

\subsubsection{Cycle Transition Matrix of Constellation Orbit }
The replenishment cycle can be modeled as $(r,q,\tau)$ policy introduced in the Section \ref{sec:inventory_model}.
Within the review period  $\tau_\text{c}$, the number of satellites will drop according to the $P_f$ in Eq.~\eqref{eq:fail_matrix_sim} $k_\text{c}$ times. Therefore, if the state distribution right after the replenishment (i.e., RAAN contact) was $\pi^{q_\text{c}}$ then $\pi^{r_\text{c}}$, which is the distribution right before the review period, becomes
\begin{equation} \label{eq:indirect_qi2ri}
    \pi^{r_\text{c}} = P_{f_\text{c}} \pi^{q_\text{c}} = \left(P_{f}\right)^{k_\text{c}} \pi^{q_\text{c}}.
\end{equation}
After having $k_\text{c}$ number of unit steps, the system undergoes the replenishment process, completing the repeated structure as follows:
\begin{equation} \label{eq:inplane_sol_eq}
\pi^{q_\text{c}} = P_{q_\text{c}}P_{f_\text{c}} \pi^{q_\text{c}},\quad \pi^{r_\text{c}} = P_{f_\text{c}}P_{q_\text{c}}\pi^{r_\text{c}}.
\end{equation}
These cycle transition matrices $P_{q_\text{c}}P_{f_\text{c}}$ and $P_{f_\text{c}}P_{q_\text{c}}$ satisfy the necessary condition for the existence of a unique stationary distribution. Therefore, any numerical method applied will efficiently compute the stationary distributions $\pi^{q_\text{c}}$ and $\pi^{r_\text{c}}$.

With this stationary distribution $\pi^{q_\text{c}}$, we can compute the average stock level during the replenishment cycle as:
\begin{equation} \label{eq:indirect_inplane_sol}
    \pi^{\text{rc}_\text{c}} = \frac{1}{k_\text{c}}\left( I + P_f + \left(P_f\right)^2 + \cdots + \left(P_f\right)^{k_\text{c}-1} \right)\pi^{q_\text{c}},
\end{equation}
which represents the average state distribution over a $\tau_\text{c}$ period. For the case of in-plane analysis, the cycle period is same as review period as
\begin{equation}
    \tau_{\text{rc}_\text{c}} = \tau_\text{c} = k_\text{c} \tau_\text{mc}.
\end{equation}

\subsubsection{Demand Distribution of Constellation Orbits to Parking Orbits}
We compute the spare demand from constellation orbits over each review period $\tau_\text{c}$, measured in units of batch size $q_\text{c}$.  Just before RAAN contact (i.e., event E), the in‐plane stock follows the distribution $\pi^{r_\text{c}}$, so the probability mass function (PMF) of the batch demand $D_\text{c}$ at the contact moment is
\begin{equation} \label{eq:inplane_demand_prob}
\begin{aligned}
    \chi_{j} &= \mathbb{P}(D_\text{c} = j | E ) \\
    &= \sum_{k=j \cdot q_\text{c}}^{\min\{(j+1) q_\text{c}-1, N_{\text{sat}_\text{c}} \}} \pi^{r_\text{c}} \left( X_\text{c} =  N_{\text{sat}_\text{c}} - k\right), \quad j=0,1,\dots
\end{aligned}
\end{equation}
Here $\chi_j$ is the probability that exactly $j$ spare batches are required when the RAANs align, and this $\chi$-vector then serves as input to the parking‐orbit analysis.  For example, if $r_\text{c} = 3$, $q_\text{c} = 2$ and ${N}_{\text{sat}_\text{c}}=5$, then
$\chi_0 = \pi^{3}(X_\text{c} = 4) + \pi^{3}(X_\text{c} = 5)$ since $D_\text{c} = 0$ when $X_\text{c} \in \left\{ 4,5\right\}$.

\subsection{Parking Orbit Analysis Method} \label{sec:event_park_analysis}
\subsubsection{Failure and Replenishment Transition Matrix of Parking Orbits}
Assuming that spare satellites do not fail, the number of parking spares decreases only when they are transferred to the constellation orbits to meet demand\footnote{If one wants to account for spare‐satellite failures, multiplying by $P_f$ at each time step yields the full transition matrix, as shown in the constellation orbit analysis.}. Using the demand distribution from Eq.~\eqref{eq:inplane_demand_prob}, the (demand‐induced) failure transition matrix of the parking orbit at RAAN contact is
\begin{equation} \label{eq:fail_matrix_park}
    P_{f_\text{p}} = \begin{bmatrix}
        \chi_{0} & 0   & \cdots & 0\\
        \chi_{1} & \chi_{0}   & \cdots & 0 \\
        \vdots & \vdots   & \ddots & \vdots \\
        1 -\sum_{i=0}^{{N}_{\text{sat}_\text{p}}-1}\chi_{i} & 1 -\sum_{i=0}^{{N}_{\text{sat}_\text{p}}-2}\chi_{i}   & \cdots & 1 \\
    \end{bmatrix},
\end{equation}
which has the same lower triangular structure as $P_{f}$ in Eq.~\eqref{eq:fail_matrix_sim}, with the failure probabilities $\nu$ replaced by the (demand-induced) failure probabilities $\chi$. In summary, letting $\pi$ be the parking state distribution before RAAN alignment, the product $P_{f_\text{p}}\pi$ gives the distribution immediately after alignment (i.e., after transferring the spares).

On the other hand, suppose a ground resupply has just been applied, and the stock level was $X_\text{p} = x$ for some $x\leq r_\text{p}$ immediately before replenishment. Then right after delivery it becomes $X_\text{p} = x+q_\text{p}$. Accordingly, we define the replenishment transition matrix $P_{q_\text{p}}$ by
\begin{equation} \label{eq:resupply_trans_matrix}
    P_{q_\text{p}} 
    =\begin{bmatrix}
    \begin{array}{c|c}
        I_{q_\text{p}} & I_{r_\text{p}+1} \\
        \mathbf{0}_{(r_\text{p}+1)\times q_\text{p}} & \mathbf{0}_{q_\text{p}\times (r_\text{p}+1)}
    \end{array}
    \end{bmatrix},
\end{equation}
and $P_{q_\text{p}} \in \mathbb{R}^{( N_{\text{sat}_\text{p}}+1) \times ( N_{\text{sat}_\text{p}}+1)}$. Note that the upper‐left block is chosen as $I_{q_\text{p}}$ to keep $P_{q_\text{p}}$ valid transition matrix. Since pre‐delivery states with $X_\text{p} > r_\text{p}$ have zero probability (enforced by $C_{r_\text{p}}^-$), any other block would give the same end result. As in Eq.~\eqref{eq:park_avail_mat}, $P_{q_\text{p}}\pi$ gives the distribution after receiving $q_\text{p}$ spares, when $\pi$ was the distribution immediately before replenishment.

Finally, to isolate probability mass above or below the reorder threshold $r_\text{p}$, define
\begin{equation} \label{eq:Cr_matrix}
    C_{r_\text{p}}^+ = \begin{bmatrix}
        I_{N_{\text{sat}_\text{p}}-r_\text{p}} & \textbf{0}_{(N_{\text{sat}_\text{p}}-r_\text{p})\times (r_\text{p}+1)} \\ \textbf{0}_{(r_\text{p}+1) \times (N_{\text{sat}_\text{p}}-r_\text{p})} & \textbf{0}_{r_\text{p}+1}
    \end{bmatrix},\ 
    C_{r_\text{p}}^- = I_{N_{\text{sat}_\text{p}}} - C_{r_\text{p}}^+
\end{equation}
Then $C_{r_\text{p}}^+\pi$ gives the distribution for $X>r_\text{p}$ and $C_{r_\text{p}}^-\pi$ the distribution for $X\leq r_\text{p}$.  These projections are key for deriving the failure and replenishment transition matrix.

\subsubsection{Transition Matrix from Delivery to Reorder} 
This section derives the transition matrix from $\pi^{q_\text{p}}$ to $\pi^{r_\text{p}}$. Within each review cycle of length $k_\text{p}$, a replenishment can arrive at any time step $i = 1, \dots, k_\text{p}$. Let $\pi^{q_\text{p}}_i$ denote the state distribution immediately after a delivery at step $i$. The average post-delivery distribution is given by
\begin{equation}
    \pi^{q_\text{p}} = \sum_{i=1}^{k_\text{p}} \eta_i \pi^{q_\text{p}}_i,
\end{equation}
where $\eta_i$ is the probability of delivery at time step $i$, whose expression is shown later in Eq.~\eqref{eq:eta_i}.

\begin{figure}[!ht]
    \centering
    \includegraphics[width=.65\textwidth]{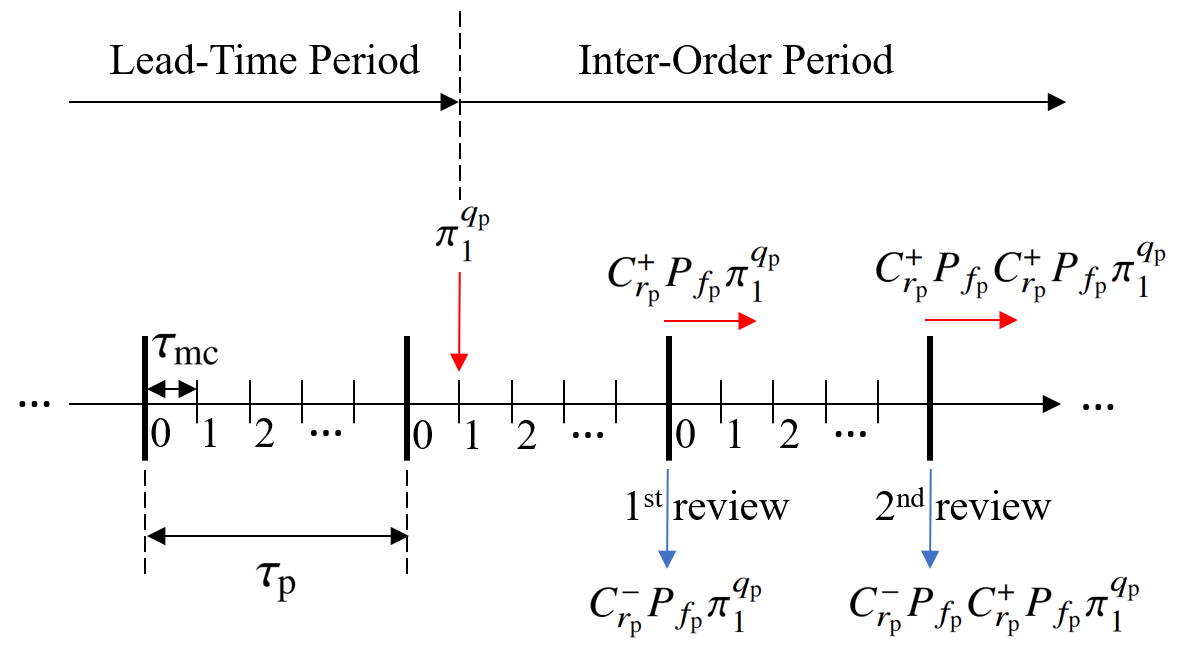}
    \caption{Transition diagram during IO period when $\pi^{q_\text{p}}_1$ happens}
    \label{fig:transition_path}
\end{figure}

After a delivery at step $i$, the state distribution $\pi^{q_\text{p}}_i$ remains unchanged until the next review at step $k_\text{p}$, since no demands occur in between. At the review, part of the population $C_{r_\text{p}}^- P_{f_\text{p}} \pi^{q_\text{p}}_i$ triggers a reorder, while the rest $C_{r_\text{p}}^+ P_{f_\text{p}} \pi^{q_\text{p}}_i$ does not. After this review, the distribution again remains fixed for the next $k_\text{p}$ steps, until the following review. At the second review, the non-reordering portion may reorder with distribution $C_{r_\text{p}}^- P_{f_\text{p}} C_{r_\text{p}}^+ P_{f_\text{p}} \pi^{q_\text{p}}_i$, or continue without reordering with distribution $C_{r_\text{p}}^+ P_{f_\text{p}} C_{r_\text{p}}^+ P_{f_\text{p}} \pi^{q_\text{p}}_i$ (see Fig.~\ref{fig:transition_path}). Enumerating all such paths to the next reorder gives
\begin{equation} \label{eq:piq_to_pir}
\begin{aligned}
    \pi^{r_\text{p}} &= \sum_{i=1}^{k_\text{p}} \sum_{j=0}^\infty C_{r_\text{p}}^- P_{f_\text{p}} \left( C_{r_\text{p}}^+ P_{f_\text{p}} \right)^{j} \eta_i \pi_i^{q_\text{p}} \\
    & = C_{r_\text{p}}^- P_{f_\text{p}} \left( I - C_{r_\text{p}}^+ P_{f_\text{p}} \right)^{-1} \sum_{i=1}^{k_\text{p}} \eta_i \pi_i^{q_\text{p}} \\
    &= C_{r_\text{p}}^- P_{f_\text{p}} \left( I - C_{r_\text{p}}^+ P_{f_\text{p}} \right)^{-1} \pi^{q_\text{p}}
\end{aligned}
\end{equation}

\subsubsection{Transition Matrix from Reorder to Delivery}
This subsection derives the transition matrix from $\pi^{r_\text{p}}$ to $\pi^{q_\text{p}}$. Under the $(r,q,\tau)$ policy, a reorder is placed only at the review point, which occurs at step 0. The probability that the replenishment arrives at step $i$ of the $j^\text{th}$ review period is given by $\rho_{i + j k_\text{p}}$, where $\rho_l$ is the probability that the lead time equals $l$ time steps.

To compute the marginal probability that a delivery occurs at step $i$ (regardless of the review cycle index $j$), we sum over all possible $j$:
\begin{equation} \label{eq:eta_i}
    \eta_i = \sum_{j=0}^\infty \rho_{i+j k_\text{p}} \text{ for } i=1,\dots,k_\text{p},
\end{equation}
Given this, the distribution after a delivery at step $i$, defined as $\pi^{q_\text{p}}_i$, can be related to $\pi^{r_\text{p}}$ by:
\begin{equation} \label{eq:eta_piq}
    \eta_i \pi^{q_\text{p}}_i = \sum_{j=0}^\infty \rho_{i+j k_\text{p}} P_{q_\text{p}}
    \left(P_{f_\text{p}}\right)^j \pi^{r_\text{p}}.
\end{equation}
where $\left(P_{f_\text{p}}\right)^j \pi^{r_\text{p}} $ gives the distribution after $j$ reviews without delivery, and $P_{q_\text{p}}$ applies the delivery transition when it occurs at step $i + jk_\text{p}$. For instance, deliveries that occur in the second review period ($j=2$) contribute $P_{q_\text{p}}\left(P_{f_\text{p}}\right)^2 \pi^{r_\text{p}}$ with probability $\rho_{i + 2k_\text{p}}$. Finally, the average post-delivery distribution is obtained by:
\begin{equation} \label{eq:pi_qp}
\begin{aligned}
    \pi^{q_\text{p}} &= \sum_{i=1}^{k_\text{p}} \eta_i \pi^{q_\text{p}}_i \\
   &= \sum_{i=1}^{k_\text{p}}\sum_{j=0}^\infty \rho_{i+j k_\text{p}} P_{q_\text{p}}
   \left(P_{f_\text{p}}\right)^j \pi^{r_\text{p}} \\
   &= \sum_{j=0}^\infty \sum_{i=0}^{k_\text{p}-1} \rho_{i+j k_\text{p}+1} P_{q_\text{p}}
   \left(P_{f_\text{p}}\right)^j \pi^{r_\text{p}}.
\end{aligned}
\end{equation}
Rearranging Eqs.~\eqref{eq:piq_to_pir} and \eqref{eq:pi_qp} yields the cycle transition matrix, from which $\pi^{q_\text{p}}$ and $\pi^{r_\text{p}}$ are obtained. For a general lead-time model, an approximated summation must typically be used. However, for the assumed lead-time distribution in Eq.~\eqref{eq:rho_lv}, an analytical expression for Eq.~\eqref{eq:pi_qp} can be derived, and the result is given in Eq.~\eqref{eq:pir_to_piq} in Appendix A.

\subsubsection{Distribution during Inter-Order (IO) Period}
This subsection derives the expression for $\pi^{\text{io}_\text{p}}$ in terms of $\pi^{q_\text{p}}$. Using the previously computed $\pi^{q_\text{p}}_i$, the distribution that avoids reordering for the next $j$ consecutive review points is given by $\left( C_{r_\text{p}}^+  P_{f_\text{p}}\right)^j \pi^{q_\text{p}}_i$ (see Fig.~\ref{fig:transition_path}). Since the state remains unchanged between reviews, the average state distribution during the IO period contributed by $\pi^{q_\text{p}}_i$ is proportional to:
\begin{equation} \label{eq:pi_io_p_1}
\begin{aligned}
    \pi^{\text{io}_\text{p}}_i &\propto \left( (k_\text{p}- i) I + \sum_{j=1}^{\infty} k_\text{p} \left(  C_{r_\text{p}}^+ P_{f_\text{p}}\right)^j \right)\pi^{q_\text{p}}_i \\
    &= \left( (k_\text{p}- i) I +  k_\text{p} C_{r_\text{p}}^+ P_{f_\text{p}} \left( I - C_{r_\text{p}}^+ P_{f_\text{p}}  \right)^{-1} \right)\pi^{q_\text{p}}_i 
\end{aligned}
\end{equation}
the first term $(k_\text{p} - i)$ reflects the remaining time steps before the first review epoch following delivery at step $i$. The second term accounts for the cumulative contribution from future review periods, each consisting of $k_\text{p}$ steps.

Aggregating over all possible delivery steps $i$, weighted by their probabilities $\eta_i$, yields the final result as
\begin{equation} \label{eq:pi_io_p_2}
\begin{aligned}
    \pi^{\text{io}_\text{p}} &= \sum_{i=1}^{k_\text{p}} \eta_i \pi^{\text{io}_\text{p}}_i \\
    &= \frac{1}{k_{\text{io}_\text{p}}} 
    \left(
    \sum_{i=1}^{k_\text{p}} (k_\text{p} - i)\eta_i \pi^{q_\text{p}}_i 
    + k_\text{p} C_{r_\text{p}}^+P_{f_\text{p}}(I - C_{r_\text{p}}^+P_{f_\text{p}})^{-1} \pi^{q_\text{p}}
    \right)
\end{aligned}
\end{equation}
where \( k_{\text{io}_\text{p}} \) is the normalization constant that ensures a valid state distribution. It also represents the expected IO period length in units of \( \tau_\text{mc} \), giving the time-converted interval as $\tau_{\text{io}_\text{p}} = k_{\text{io}_\text{p}} \tau_\text{mc}$.

\subsubsection{Distribution during Lead-Time Period} 
This subsection derives the expression for $\pi^{\text{lt}_\text{p}}$ in terms of $\pi^{r_\text{p}}$. The probability that the replenishment has not arrived by the $\ell^\text{th}$ time step is 
\begin{equation}
    \rho^c_\ell = 1 - \sum_{i=1}^\ell \rho_i, \quad \ell = 0,1,\dots,
\end{equation}
with $\rho_0^c = 1$. Therefore, the weighted distribution that is still awaiting delivery at the $i^\text{th}$ time step of the $j^\text{th}$ review period after the reorder moment is $ \rho_{i+jk_\text{p}}^c \left(P_{f_p}\right)^j \pi^{r_\text{p}}$. Then, the average distribution during the LT period becomes
\begin{equation} \label{eq:pi_lt_p1}
\begin{aligned}
    \pi^{\text{lt}_\text{p}} 
    &= \frac{1}{k_{\text{lt}_\text{p}}}  \sum_{j=0}^\infty \sum_{i=0}^{k_\text{p}-1} \rho_{i+jk_\text{p}}^c \left(P_{f_p}\right)^j \pi^{r_\text{p}} \\
    &= \frac{1}{k_{\text{lt}_\text{p}}} \left[ \left(\sum_{i=0}^{k_\text{p}-1} \rho_i^c\right) I + \left(\sum_{i=0}^{k_\text{p}-1} \rho_{i+k_\text{p}}^c \right)P_{f_\text{p}} +\cdots \right] \pi^{r_\text{p}}
\end{aligned}
\end{equation}
where \( k_{\text{lt}_\text{p}} \) is the normalization constant and also represents the expected length of the LT period in time steps. As before, we can derive the analytical expression of Eq.\eqref{eq:pi_lt_p1} based on the assumed lead-time distribution in Eq.\eqref{eq:rho_lv}, and the resulting expression is given in Eq.~\eqref{eq:pi_lt_p2} in Appendix A.

\subsubsection{Distribution during Every Replenishment Cycle} \label{sec:park_rc}
Finally, since we have derived the state distributions during the IO and LT periods, along with their respective durations, the average state distribution in the parking orbit over a complete replenishment cycle under the indirect resupply policy is given by:
\begin{equation} \label{eq:indirect_park_sol}
\pi^{\text{rc}_\text{p}} = \frac{k_{\text{io}_\text{p}}}{k_{\text{io}_\text{p}} + k_{\text{lt}_\text{p}}}\pi^{\text{io}_\text{p}} + \frac{k_{\text{lt}_\text{p}}}{k_{\text{io}_\text{p}} + k_{\text{lt}_\text{p}}}\pi^{\text{lt}_\text{p}},
\end{equation}
and the corresponding average duration of one replenishment cycle in the parking orbit is:
\begin{equation}
    \tau_{\text{rc}_\text{p}} = \tau_{\text{io}_\text{p}} + \tau_{\text{lt}_\text{p}}.
\end{equation}

\subsubsection{Parking Spares Availability Distribution} \label{sec:park_avail}
From the parking-orbit analysis, we compute the availability distribution at the RAAN-contact event $E$ (or more specifically, just before parking spares can be distributed to the in-plane orbit). As stated in Eq.\eqref{eq:parking_avail}, we require the parking-spares distribution at the moment $E$ occurs, which may fall within either the IO period or the LT period. Denoting the corresponding distribution for each case by $\pi^{\text{io}\text{p}|E}$ and $\pi^{\text{lt}\text{p}|E}$, whose closed-form expressions are derived in Eqs.\eqref{eq:pi_iop_E} and \eqref{eq:pi_ltp_E} in the Appendix A, the average distribution conditioned on the event $E$ during the entire period becomes:
\begin{equation} \label{eq:pi_rc_E}
    \pi^{\text{rc}_\text{p}|E} = 
    \frac{k_{\text{io}_\text{p}|E}}{k_{\text{io}_\text{p}|E} + k_{\text{lt}_\text{p}|E}} \pi^{\text{io}_\text{p}|E} +
    \frac{k_{\text{lt}_\text{p}|E}}{k_{\text{io}_\text{p}|E} + k_{\text{lt}_\text{p}|E}}
    \pi^{\text{lt}_\text{p}|E}
\end{equation}
where $k_{\text{io}_\text{p}|E}$ and $k_{\text{lt}_\text{p}|E}$ are the normalization constants of each distribution. Then, the parking availability probability in Eq.~\eqref{eq:parking_avail} can be computed as:
\begin{equation} \label{eq:kappa_eq}
    \kappa_j = \sum_{k=j}^{N_{\text{sat}_\text{p}}}\pi^{\text{rc}_\text{p}|E} (X_\text{p} = k)
\end{equation}

\subsection{Flow of Indirect Strategy Analysis}
We split the coupled in-plane and parking-orbit analysis into two parts, as explained in the preceding subsections. Solving them independently yields inconsistent results, since they are coupled through $\chi$ and $\kappa$, so we use fixed-point iteration, starting from 100 \% parking availability($\kappa = \mathbbm{1}$), to reach a consistent solution.

Although the analysis involves multiple steps, it takes only a few milliseconds on a standard computer. The detailed procedure for the indirect-resupply method is summarized in Table \ref{alg:cap}.

\begin{algorithm}
\caption{Fixed Point Iteration for the Indirect Strategy}\label{alg:cap}
\begin{algorithmic}
\Require Constellation Configuration, Probability Model
\State $\tau_\text{c} \gets$ Eq.~\eqref{eq:t_plane_indirect} and $\tau_\text{p} \gets$ Eq.~\eqref{eq:t_park_indirect}
\State $P_{f} \gets$ Eq.~\eqref{eq:fail_matrix_sim} and $P_{f_\text{c}} \gets$ Eq.~\eqref{eq:indirect_qi2ri} 
\State $P_{q_\text{p}} \gets$ Eq.~\eqref{eq:resupply_trans_matrix} and $C_{r_\text{p}}^-, C_{r_\text{p}}^+ \gets$ Eq.~\eqref{eq:Cr_matrix}
\State $\kappa \gets \mathbbm{1}$ and $k\gets 0$
\While{ $ |\kappa^{k+1} - \kappa^{k}| >\varepsilon $ or $k < k^{\max}$}
    \State $P_{q_\text{c}}\gets$ Eq.~\eqref{eq:park_avail_mat}
    \State $\pi^{q_\text{c}}$, $\pi^{r_\text{c}} \gets$ Eq.~\eqref{eq:inplane_sol_eq}
    and $\pi^{\text{rc}_\text{c}} \gets$ Eq.~\eqref{eq:indirect_inplane_sol}
    \State $\chi \gets$ Eq.~\eqref{eq:inplane_demand_prob} and $P_{f_\text{p}} \gets $ Eq.~\eqref{eq:fail_matrix_park}
    \State $\pi^{q_\text{p}}$, $\pi^{r_\text{p}} \gets$ Eqs.~\eqref{eq:piq_to_pir},~\eqref{eq:pi_qp}
    \State $\pi^{\text{io}_\text{p}},\pi^{\text{lt}_\text{p}} \gets$ Eqs.~\eqref{eq:pi_io_p_2},~\eqref{eq:pi_lt_p1} and $\pi^{\text{rc}_\text{p}} \gets$  Eq.~\eqref{eq:indirect_park_sol}
    \State $\pi^{\text{io}_\text{p}|E},\pi^{\text{lt}_\text{p}|E} \gets$ Eqs.~\eqref{eq:pi_iop_E},~\eqref{eq:pi_ltp_E} and 
    $\pi^{\text{rc}_\text{p}|E} \gets$ Eq.~\eqref{eq:pi_rc_E} 
    \State $\kappa \gets$ Eq.~\eqref{eq:kappa_eq} and $k \gets k + 1$
\EndWhile
\end{algorithmic}
\end{algorithm}

\section{Performance Evaluation of Spare Management Policy} \label{sec4}
With the stationary solution \(\pi^{\mathrm{rc}_{(\cdot)}}\) in hand, we can evaluate general performance metrics. The two most common metrics are operational cost and resilience, which trade off with each other.

\subsection{Cost Model of Indirect Strategy}
The total expected operating cost per unit time $C_\text{total}$ is
\begin{equation}
    C_\text{total} = C_\text{build} + C_\text{hold} + C_\text{trans} + C_\text{launch} 
\end{equation}
where $C_\text{build}$ is the expected manufacturing cost of spares per unit time,
$C_\text{hold}$ is the expected holding cost per unit time, 
$C_\text{trans}$ is the expected orbital maneuver cost per unit time, and
$C_\text{launch}$ is the expected launch cost per unit time.
 
Firstly, $C_\text{build}$ is defined as:
\begin{equation}
    C_\text{build} = \frac{1}{\tau_{\text{rc}_\text{p}}}c_\text{build}  N_{\text{orbit}_\text{p}} \cdot q_\text{c} \cdot q_\text{p} 
\end{equation}
where $c_\text{build}$ is the manufacturing cost of a spare satellite. This reflects that $q_\text{c}\cdot q_\text{p}$ spares ($=q_\text{p}$ batches) are launched to each of the $N_{\text{orbits}_\text{p}}$ parking orbits every $\tau_{\text{rc}_\text{p}}$

The holding cost represents the penalty for maintaining an excessive number of spares in orbit. It accounts for station-keeping costs, depreciation, and failure risk. It can be modeled as:
\begin{equation}
    C_\text{hold} = c_{\text{hold}_\text{c}}  N_{\text{orbit}_\text{c}} \sum_{i=\bar{N}_\text{sat}+1}^{{N}_{\text{sat}_\text{c}}} \left( i - \bar{N}_\text{sat} \right) 
    \pi^{\text{rc}_\text{c}}(X_\text{c} = i)
    + c_{\text{hold}_\text{p}} N_{\text{orbit}_\text{p}} q_\text{c}\sum_{i=0}^{{N}_{\text{sat}_\text{p}}} i 
    \pi^{\text{rc}_\text{p}}(X_\text{p} = i)
\end{equation}
where $c_{\text{hold}_\text{c}}$ and $c_{\text{hold}_\text{p}}$ are the holding cost of a spare satellite per unit time for constellation orbits and parking orbits, respectively.

The orbit transfer cost is modeled as follows:
\begin{equation}
    C_\text{trans} = \frac{N_{\text{orbit}_\text{c}}}{\tau_{\text{rc}_\text{c}} q_\text{c}} (c_\text{fuel}\cdot m_\text{fuel} + c_\text{trans})  \sum_{i=0}^{{N}_{\text{sat}_\text{c}}} i \left( \pi^{q_\text{c}}(X_\text{c}=i) -  \pi^{r_\text{c}}(X_\text{c}=i)\right)
\end{equation}
where \(c_{\rm fuel}\) is the cost of fuel per unit mass, \(c_{\rm trans}\) covers non‐fuel transfer costs (e.g. bus build, transfer risk), and the fuel mass \(m_{\rm fuel}\) is given by Eqs.~\eqref{eq:delta_v} and \eqref{eq:fuel_mass}, with \(m_{\rm dry}=q_\text{c}\,m_{\rm sat}+m_{\rm bus}\). 

Lastly, the expected launch cost $C_\text{launch}$ is modeled under two scenarios, depending on whether rideshare opportunities for parking orbits are available:
\begin{equation}
    C_\text{launch} = 
    \begin{cases} 
    \dfrac{N_{\text{orbit}_\text{p}}}{\tau_{\text{rc}_\text{p}}} \min\left\{  c_\text{lv,unit} m_\text{total}   ,\ c_\text{lv,full}
    \right\}, & \text{if rideshare available}, \\
    \dfrac{N_{\text{orbit}_\text{p}}}{\tau_{\text{rc}_\text{p}}} c_\text{lv,full}, & \text{if rideshare unavailable}. 
\end{cases}
\end{equation}
where $c_\text{lv,unit}$ is the launch cost to LEO per unit mass, 
$c_\text{lv,full}$ is the discounted price for reserving the full payload capacity of the launch vehicle, and 
$m_\text{total} = (m_\text{fuel} + m_\text{dry}) \cdot q_\text{p}$ is the total mass of the $q_\text{p}$ batches. 
The minimum operator in the first case reflects the economic choice between paying the per-kilogram rate 
and purchasing a full-vehicle contract when rideshare missions to the target orbit are available. 
The second case assumes that rideshare is not offered for the desired orbit, so it must always purchase the full vehicle.

\subsection{Resilience Model of Indirect Resupply Strategy}
The proper resilience metric should capture both agility (how quickly the system recovers to nominal capacity) and robustness (the depth of capability loss when operating below nominal) \cite{dod_resilience_2011, Ron_resilience}. In practice, this is often measured by the area under the nominal‐capacity line that is lost over time, i.e.\ the time‐integral of the deficit. In our discrete‐time Markov model this becomes the expected shortage $S_\text{c}$
\begin{equation}
    S_\text{c} = \sum_{i=0}^{\bar{N}_\text{sat}} (\bar{N}_{\text{sat}} - i)\cdot \pi^{\text{rc}_\text{c}} (X_\text{c}=i),
\end{equation}
which multiplies each deficit $\bar{N}_{\text{sat}} - i$ by the fraction of time spent in state $i$ as $\pi^{\text{rc}_\text{c}} (X_\text{c}=i)$.

In a parking orbit there is no nominal capacity level, so the expected-shortage metric does not apply. Instead, we define resilience as the out-of-stock probability:
\begin{equation} \label{eq:res_park}
    \mathbb{P}(X_\text{p} = 0)=  \pi^{\text{rc}_\text{p}} (X_\text{p}=0) 
\end{equation}
which is the long-run fraction of time the parking orbit has zero spares.

\subsection{Optimization Problem for Indirect Resupply Strategy}
There are multiple ways to formulate the optimization problem, but here we focus on minimizing the total operating cost of the spare policy while enforcing resilience and launch‐vehicle constraints. The problem can be stated as
\begin{equation} \label{eq:opt_formulation}
\begin{aligned}
    \min_{x}\quad & C_\text{total}\\
    \text{s.t.}\quad & g_1 = S_\text{c} - \varepsilon_1 \leq 0 \\
    \quad & g_2 = \mathbb{P}(X_\text{p} = 0) - \varepsilon_2 \leq 0 \\
    \quad & g_3 = m_\text{total} - m_\text{payload} \leq 0 \\
    \quad & q_\text{c}, r_\text{c}, q_\text{p}, r_\text{p}, N_{\text{orbit}_\text{p}} \in \mathbb{Z}^+ \\
\end{aligned}
\end{equation}
where $x = (q_\text{c}, r_\text{c}, q_\text{p}, r_\text{p}, N_{\text{orbit}_\text{p}}, h_\text{p})$, 
and $\varepsilon_1$ and $\varepsilon_2$ are user-defined thresholds. 
The $g_3$ constraint ensures that the total mass does not exceed the vehicle’s maximum payload capacity $m_\text{payload}$.  
By solving this optimization problem, one can estimate the approximate cost required 
to maintain the satellite constellation.

\section{Numerical Validation of the Analysis Method} \label{sec5}
\subsection{Numerical Validation Set-up}
The analytical model developed above enables efficient evaluation of spare policies even for mega‐scale constellations, but it relies on simplifying assumptions such as independent and identically distributed stock levels in both constellation and parking orbits. To validate the proposed method, we follow the approach introduced in \cite{jakob2019optimal}.

When validating the model, we could directly compare histograms of simulated stock levels with \(\pi^{\mathrm{rc}_{(\cdot)}}\), but summarizing those differences in a single metric is difficult. Instead, we compare the mean stock levels of in‐plane and parking orbits, the expected shortage of the constellation orbit, and the out‐of‐stock probability of the parking orbit, which are key metrics for performance evaluation. The mean stock level is computed as
\begin{equation}
    M_{(\cdot)} = \sum_{i=0}^{N_{\text{sat}_{(\cdot)}}} i\cdot \pi^{\text{rc}_{(\cdot)}} (X_{(\cdot)}=i).
\end{equation}

Finally, we construct 100 unique test cases using Latin hypercube sampling over the parameters in Table~\ref{tab:trade_space_lhs}, with fixed parameters given in Table~\ref{tab:fixed_sim_para}. For each test case, we run a 20-year simulation 100 times and report the averaged statistics.

\subsection{Numerical Validation Results} \label{sec:num_val}
As noted in our earlier work \cite{han2024analysis}, the proposed analysis can lose accuracy when the i.i.d. assumption for parking orbits breaks down. This occurs if parking spares are out of stock for extended periods, which typically happens when the satellite failure rate exceeds the replenishment capacity and disrupts the conventional saw-tooth stock profile. Consequently, the heuristic threshold $\mathbb{P}(X_\text{p} = 0) < 1/(N_{\text{sat}_\text{p}}+1)$, derived from an even-distribution assumption, is used here to select valid cases for comparison.

\begin{figure}[!ht]
    \centering
    \includegraphics[width=.45\textwidth]{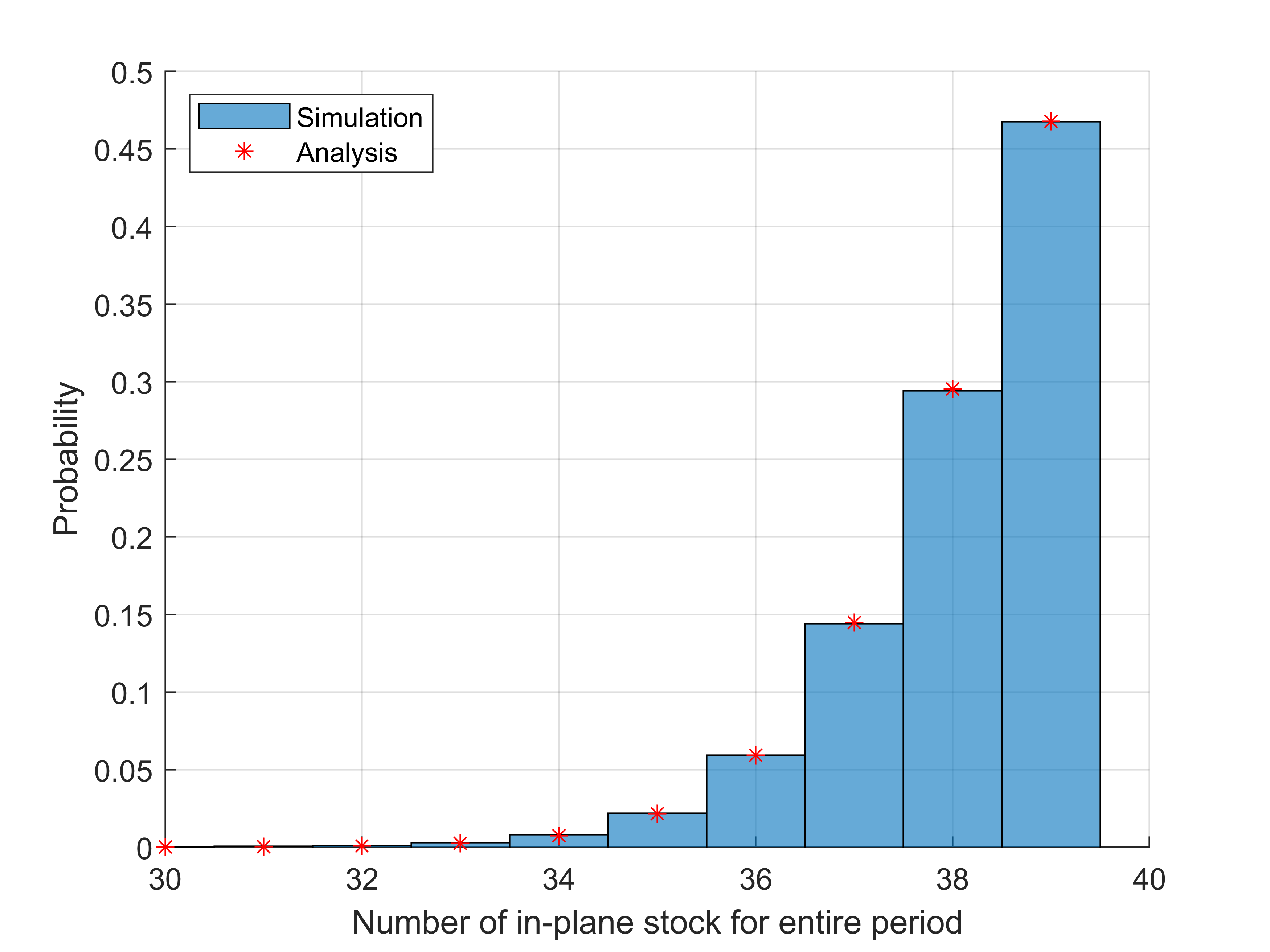}
    \caption{Comparison of $\pi^{\text{rc}_\text{c}}$ for a representative case near the maximum error level.}
    \label{fig:worst_case_result_c}
\end{figure}

\begin{figure}[!ht]
    \centering
    \includegraphics[width=.45\textwidth]{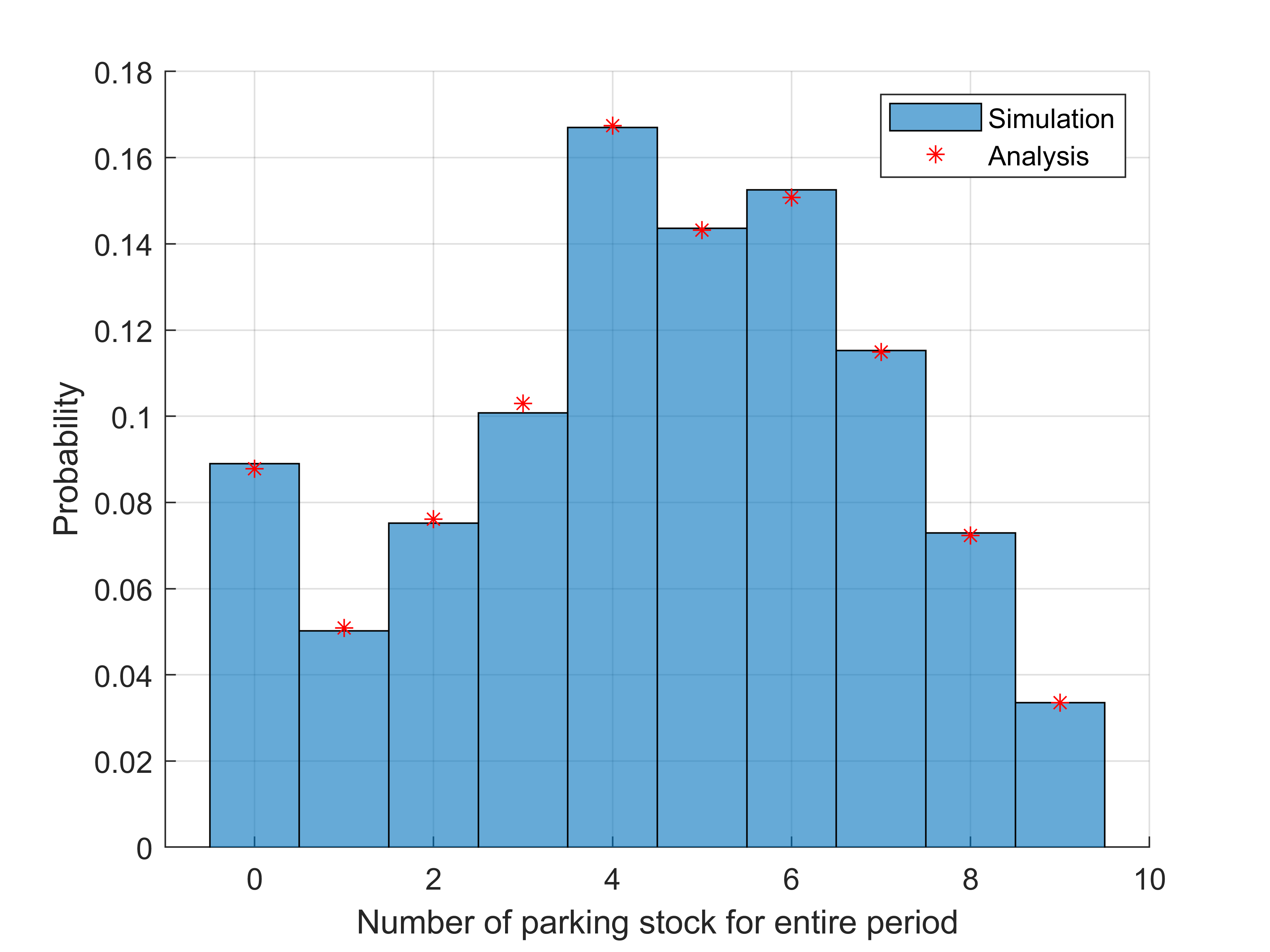}
    \caption{Comparison of $\pi^{\text{rc}_\text{p}}$ for a representative case near the maximum error level.}
    \label{fig:worst_case_result_p}
\end{figure}

The test results are summarized in Table~\ref{tab:error_list}. For error computation, relative error is used for $M_\text{c}$, $M_\text{p}$, and $S_\text{c}$, while absolute error is used for $\mathbb{P}(X_\text{p} = 0)$. The fidelity of the analytical framework is confirmed by the results, which show a 95th percentile error of less than 1\% across all performance metrics. The minor deviations between the analytical and simulated results are attributed to the model's i.i.d. assumption for the parking orbits and the statistical noise inherent in finite-duration Monte Carlo simulations. The higher relative errors reported for $M_\text{p}$ and $S_\text{c}$ compared to $M_\text{c}$ are primarily an artifact of their small nominal values, which is a known characteristic of the relative error metric. The overall validity of the method is further demonstrated in Figs.~\ref{fig:worst_case_result_c} and~\ref{fig:worst_case_result_p}, which show a close alignment of the constellation- and parking-orbit state distributions for a representative near-worst-case test case.

In terms of computational cost, the full Monte Carlo simulation required several hours to complete each test case, whereas the proposed method produced a solution in under a second. This demonstrates that the proposed method is both accurate and efficient, making it suitable for use as an inner loop in the optimization process.

\begin{table}[hbt!]
\centering
\caption{Fixed simulation parameters}
\label{tab:fixed_sim_para}
\begin{tabular}{lc c c c}
\hline \hline
Parameter & Notation & Value & Unit \\
\hline

Markov time step
  & $\tau_\text{mc}$ 
  & $0.5$ 
  & days \\

Constellation orbit altitude ($a_\text{c} - R_\oplus$)
  & $h_\text{c}$ 
  & $1200$ 
  & km \\
  
Inclination of orbit planes
  & $i$ 
  & $50$ 
  & deg \\

Number of constellation orbits
  & $N_{\text{orbit}_\text{c}}$ 
  & $40$ 
  & orbits \\

Nominal satellites per plane
  & $\bar{N}_\text{sat}$ 
  & $40$ 
  & satellites \\

\hline
\hline
\end{tabular}
\end{table}

\begin{table}[hbt!]
\centering
\caption{Bound of sampled simulation parameters}
\label{tab:trade_space_lhs}
\begin{tabular}{lc c c c}
\hline \hline
Parameter & Notation & Bounds & Unit \\
\hline

Annual satellite failure rate
  & $\lambda_{\text{sat,yr}}$ 
  & $\left[ 0.001,\  0.5\right]$ 
  & failures/satellite/year \\

Launch order processing time 
  & $\tau_{\text{lv}}$ 
  & $\left[ 0,\  60\right]$ 
  & days \\

Mean exponential launch lead time
  & $\mu_{\text{lv}}$ 
  & $\left[ 5,\  60\right]$ 
  & days \\

Order size for in-plane spares
  & $q_{\text{c}}$ 
  & $\left[ 1,\  20\right]$ 
  & satellites \\

Order size for parking spare batch
  & $q_{\text{p}}$ 
  & $\left[ 1,\  20\right]$ 
  & batches \\

Reorder point for in-plane spares
  & $r_{\text{c}}$ 
  & $\left[\bar{N}_\text{sat}- 5,\  \bar{N}_\text{sat}+5\right]$ 
  & batches \\  

Reorder point for parking spare batch
  & $r_{\text{p}}$ 
  & $\left[0,\  10\right]$ 
  & batches \\  

Parking orbit altitude ($a_\text{p} -R_\oplus$)
  & $h_\text{p}$ 
  & $\left[ 500,\  1100\right]$ 
  & km \\

Number of parking orbits
  & $N_{\text{orbit}_\text{p}}$ 
  & $\left[ 1,\  20\right]$ 
  & orbits \\

\hline
\hline
\end{tabular}
\end{table}

\begin{table}[hbt!]
\centering
\caption{Error between proposed method and simulation results}
\label{tab:error_list}
\begin{tabular}{lc c c }
\hline \hline
Parameter & Mean & P95\\
\hline
Relative error of $M_\text{c}$ 
  & $0.012 \ \%$ 
  & $0.035 \ \%$ \\

Relative error of $M_\text{p}$ 
  & $0.172 \ \%$ 
  & $0.432 \ \%$ \\
  
Relative error of $S_\text{c}$ 
  & $0.191 \ \%$ 
  & $0.794 \ \%$ \\

Absolute error of $\mathbb{P}(X_\text{p} = 0)$ 
  & $0.006 \ \text{p.p}$  
  & $0.019 \ \text{p.p}$  \\
\hline
\hline
\end{tabular}
\end{table}

\subsection{Comparison with Existing Method}
\begin{figure}[!ht]
    \centering
    \includegraphics[width=.45\textwidth]{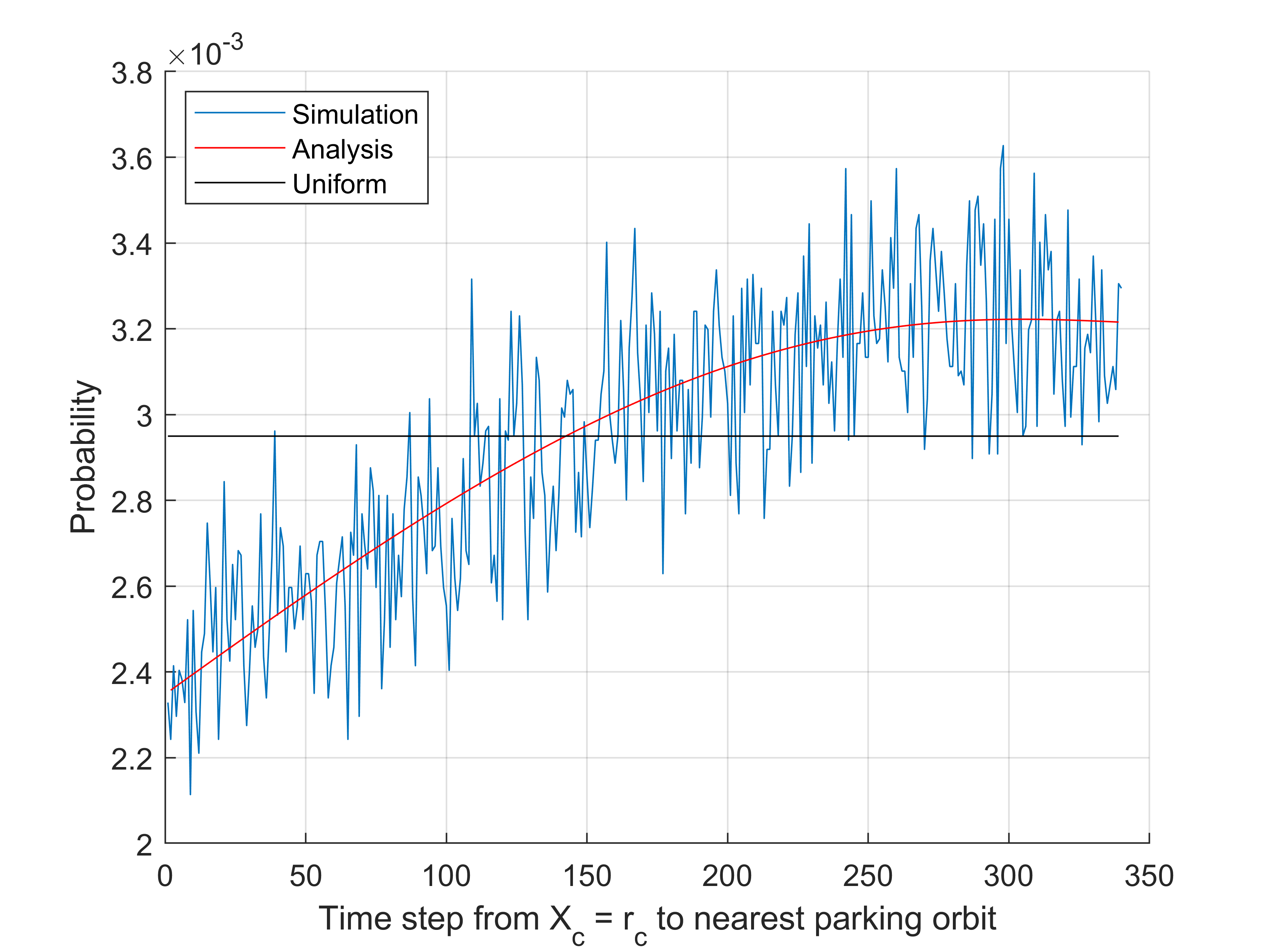}
    \caption{Comparison of lead-time distributions from simulation, the proposed analysis, and the uniform model}
    \label{fig:lead_time_test}
\end{figure}

\begin{figure}[!ht]
    \centering
    \includegraphics[width=.45\textwidth]{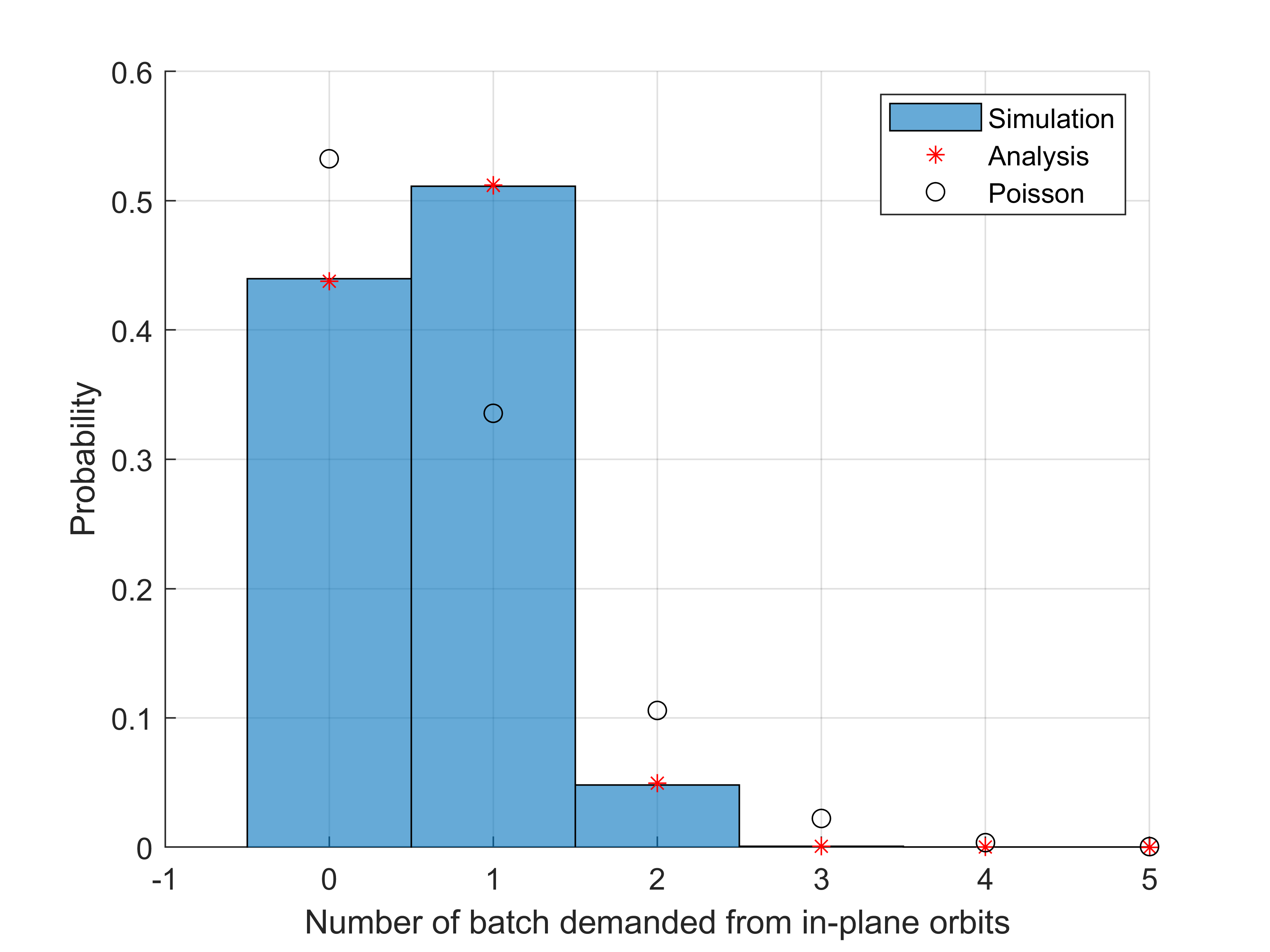}
    \caption{Histogram of in-plane spares demand at RAAN alignment for a representative case}
    \label{fig:demand_test}
\end{figure}

The previous work \cite{jakob2019optimal} assumes that the demand for parking spares follows a Poisson process and that the lead time from parking orbits to constellation orbits follows a uniform distribution, based on the assumption of a high order fill rate. Specifically, the lead time is defined as the duration from the moment $X_\text{c} \leq r_\text{c}$ to the next RAAN alignment with the nearest parking orbit, and the following parameter is used for Poisson demand \cite{jakob2019optimal}:
\begin{equation}
    \lambda_\text{demand} = \frac{N_{\text{orbit}_\text{c}} \bar{N}_{\text{sat}}}{q_\text{c} N_{\text{orbit}_\text{p}}} \lambda_\text{sat} \tau_\text{c}.
\end{equation}
However, our simulations reveal significant deviations from these assumptions. Figure~\ref{fig:lead_time_test} compares the lead-time distribution from simulation (blue line), the proposed analysis (red line), and the uniform model (black line). The simulated distribution is skewed and clearly non-uniform, which is not captured by the uniform model. Similarly, Fig.~\ref{fig:demand_test} demonstrates that the Poisson distribution (black circles) is a poor fit for the simulated spare demand.

These discrepancies persist even in near-ideal scenarios with low failure rates and high parking-spare availability, indicating a fundamental limitation in the prior model's ability to represent the system's stochastic interplay. In contrast, our proposed analysis (red line and circles) shows excellent agreement with the simulation data in both figures.

\section{Optimization of Spare Management Policy} \label{sec6}
One key application of the proposed analysis method is design optimization. Unless otherwise noted, the constellation configuration and baseline stochastic parameters in this section are the same as those in Sec.~\ref{sec5}; the corresponding values are given in Table~\ref{tab:fixed_sim_para}, and the additional cost parameters used for optimization are listed in Table~\ref{tab:opt_para}.  An optimization problem can be formulated to guide early-phase spare strategy decisions, such as monthly replenishment needs or cost estimates. In this section, we optimize the indirect strategy using the formulation in Eq.~\eqref{eq:opt_formulation} with real-world parameters. The resulting cost is then compared to that of the direct strategy. The direct-strategy results reported in this section are obtained using the corresponding $(r,q)$ direct-insertion model, which can be derived as a special case of the present framework under continuous review, as summarized in Appendix B. An alternative independent derivation and a more detailed discussion are given in \cite{han2025direct}.

For the direct strategy, we assume the use of a small launch vehicle, specifically Rocket Lab’s Electron, while the indirect strategy employs a heavy launch vehicle, namely SpaceX’s Falcon 9, with launch costs referenced from \cite{SpaceInsider_cost,spacex_rideshare}. Although the lead time for both vehicles can, in principle, be as short as two days \cite{ainvest_rocketlab_2025,spaceflightnow_falcon9_2025}, such estimates are optimistic and do not account for contract processing or shipping delays. To reflect a more realistic timeline, we adopt $(\mu_\text{lv},\tau_\text{lv}) = (20,20)$ days for Falcon 9 and $(10,10)$ days for Electron, making the Falcon 9 lead time roughly twice as long. For Electron we assume $c_\text{lv,full} = 7.5$ M\$ and $m_\text{payload} = 300$ kg.

Regular rideshare opportunities are generally limited to low altitudes with common inclinations such as sun-synchronous or ISS orbits. Since our constellation orbit does not fall into this category, rideshare is not feasible for the direct strategy, and Electron with a full-vehicle contract is assumed instead. For the indirect strategy, Falcon 9 is considered for parking-orbit replenishment, as its unit cost is substantially lower than Electron and rideshare may be feasible at lower altitudes. Accordingly, we evaluate two indirect models: one with rideshare access and one without. Under the no-rideshare assumption, using a heavy LV for the direct strategy is economically dominated because it requires a full-vehicle contract for a small replenishment need and also increases holding cost through oversized batches. The remaining parameters are summarized in Table~\ref{tab:fixed_sim_para} and Table~\ref{tab:opt_para}.

\begin{table}[hbt!]
\centering
\caption{Parameters for the optimization}
\label{tab:opt_para}
\begin{tabular}{lc c c c}
\hline \hline
Parameter & Notation & Value & Unit \\
\hline

Satellite manufacturing cost
  & $c_{\text{build}}$ 
  & $0.5$ 
  & M\$/satellite \\

In-orbit spares annual holding cost
  & $c_{\text{hold}_\text{c}}$ 
  & $0.5$ 
  & M\$/satellite/year \\

Parking spares annual holding cost
  & $c_{\text{hold}_\text{p}}$ 
  & $0.5$ 
  & M\$/satellite/year \\

Launch cost per unit mass (Falcon 9)
  & $c_{\text{lv,unit}}$ 
  & $6500$ 
  & \$/kg \\

Discounted cost for full contract (Falcon 9)
  & $c_{\text{lv,full}}$ 
  & $67$ 
  & M\$ \\

Maximum payload launch capacity (Falcon 9)
  & $m_{\text{payload}}$ 
  & $18500$ 
  & kg \\

Fuel cost per mass
  & $c_{\text{fuel}}$ 
  & $0.001$ 
  & M\$/kg \\

Non-fuel transfer cost
  & $c_{\text{trans}}$ 
  & $0.5$ 
  & M\$ \\

Mass of satellite
  & $m_\text{sat}$ 
  & $150$ 
  & kg \\

Mass of transfer bus
  & $m_\text{bus}$ 
  & $100$ 
  & kg \\

Effective exhaust velocity 
  & $v_\text{ex}$ 
  & $2.16$ 
  & km/s \\

\hline
\hline
\end{tabular}
\end{table}

\subsection{Baseline Scenario}
As a baseline scenario, we consider a moderate failure rate of $\lambda_\text{sat} = 0.05$, with resilience thresholds $\varepsilon_1 = 0.25$ for $S_c$ and $\varepsilon_2 = 1/(N_{\text{sat}_\text{p}}+1)$ for the out-of-stock probability. Optimization is performed using a genetic algorithm. For the indirect strategy, the optimal solution is $(q_c^\ast, r_c^\ast, q_p^\ast, r_p^\ast, N_{\text{orbit}_p}^\ast, h_p^\ast) = (4,40,23,2,1,735)$, regardless of rideshare availability. In this case, using the full Falcon 9 contract minimizes cost, as the savings in launch cost outweigh the increase in holding cost from parking spares. The direct strategy, which has only two integer design variables $(r,q)$, yields $(r^\ast,q^\ast)=(39,2)$.

As shown in Table~\ref{tab:base_scn}, the indirect strategy achieves a 53\% reduction in total cost, primarily due to the much lower launch cost per satellite enabled by Falcon 9. Although it incurs additional holding and transfer costs, these are offset by the launch savings. The build cost $C_\text{build}$ remains similar between the two strategies because the same resilience requirement under the same failure rate leads to a similar number of spare satellites being produced. Lastly, the direct strategy naturally avoids transfer operations, but due to the small design space, i.e., $q$ and $r$, the $S_\text{c}$ constraint is satisfied with an excessive margin.

This baseline case shows that, under realistic launch cost assumptions, the indirect strategy can achieve the required resilience at significantly lower cost than the direct strategy.

\begin{table}[hbt!]
\caption{Summary of results for representative scenarios.}
\label{tab:base_scn}
\centering
\begin{tabular}{l c c cccc c cc}
\hline \hline
   \multirow{2}{*}{Policy} & \multirow{2}{*}{$C_\text{total}$ [M\$/day]} && \multicolumn{4}{c}{Detailed Costs [M\$/day]} && \multicolumn{2}{c}{Constraints} \\
    \cline{4-7} \cline{9-10} 
 &  && $C_\text{build}$ & $C_\text{hold}$ &$ C_\text{launch}$ & $C_\text{trans}$ &&  $S_\text{c}$ &  $\mathbb{P}(X_\text{p} = 0)$ \\
\hline\hline
Direct(no-RS) &  0.9547  && 0.1094 & 0.0246 & 0.8207 & N/A && 0.0591 & N/A \\
Indirect(no-RS) &  0.4479 && 0.1082 & 0.1507 & 0.1575 & 0.0316 && 0.2387 & 0.0286 \\
Indirect(RS) &  0.4479 && 0.1082 & 0.1507 & 0.1575 & 0.0316 && 0.2387 & 0.0286 \\
\hline \hline
\end{tabular}
\vspace{4pt}\\
\footnotesize{Note: RS = rideshare available; no-RS = rideshare unavailable.}
\end{table}

\subsection{Sensitivity to failure rate}
Satellite failure rates vary over time and depend on system scale. Smaller satellites generally exhibit higher failure rates than larger ones \cite{2010_Gregory,langer2018reliability,jacklin2019small}. Since it is difficult to define a single representative value, we test a wide range of failure rates from 0.01 to 0.5 failures per year.

\begin{figure}[!ht]
    \centering
    \includegraphics[width=.45\textwidth]{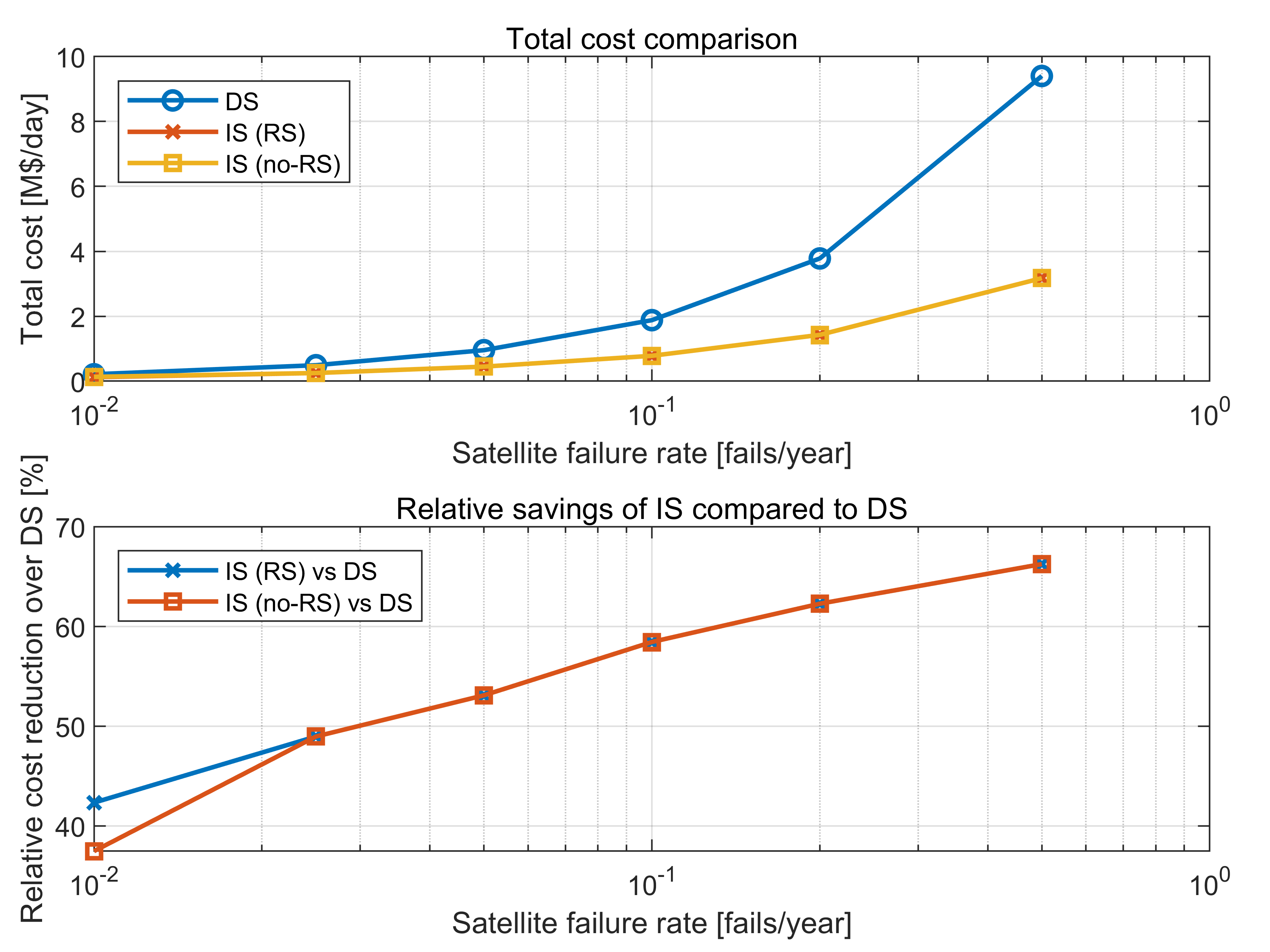}
    \caption{Relative $C_\text{total}$ savings of IS compared to the DS across different satellite failure rates.}
    \label{fig:fail_opt}
\end{figure}

In this experiment, only the satellite failure rate is varied, while all other parameters are kept the same as in the baseline scenario. Figure~\ref{fig:fail_opt} shows the total cost of each strategy and the relative savings of the indirect strategy (IS) compared to the direct strategy (DS). Higher values indicate greater cost reduction of IS over DS. IS is more cost efficient than DS across the entire failure rate range. At very low failure rates (around 0.01), relatively few spares are needed to satisfy the $S_\text{c}$ constraint. In this case, rideshare can reduce cost, while full-contract launches provide little advantage because the additional holding cost offsets the savings. As the failure rate increases, more spares and more frequent launches are required. Since launch cost dominates the budget, the advantage of IS grows with failure rate.

These results indicate that indirect replenishment via parking orbits becomes an attractive architecture for large-scale constellations with high spare demand. Accordingly, when developing actual replenishment logic for such systems, it is reasonable to consider an indirect channel alongside direct insertion. This study focuses on early-phase architecture assessment; optimizing operational logic within the selected architecture is a separate next step.

\section{Conclusion} \label{sec7}
In this paper, we developed a Markov chain-based framework for the detailed analysis and design of an indirect spare management policy in large-scale constellations. We modeled in-plane and parking orbits as $(r,q,\tau)$ systems, solved their coupled stationary distributions via fixed-point iteration, and derived expressions for cost and resilience metrics. Building on this fast, accurate analysis, we formulated an optimization problem and solved it with a genetic algorithm to minimize total operating cost under resilience constraints. Using this framework, we compared the indirect strategy to the direct strategy and characterized its behavior through sensitivity analysis. Our results show that the indirect strategy is advantageous when it can leverage lower launch costs. Finally, this framework can be extended to other constellation configurations (e.g., asymmetric and heterogeneous) and alternative replenishment policies, which we will explore in future work.

\section*{Appendix A: Closed-Form Expressions for Eqs.~\eqref{eq:pi_qp},~\eqref{eq:pi_lt_p1}, and~\eqref{eq:pi_rc_E}}
In this appendix, we derive the analytical expression using the explicit equation of the exponential lead time distribution.

\subsection*{Closed-Form Expression for Eq.~\eqref{eq:pi_qp}}
\begin{figure}[!ht]
    \centering
    \includegraphics[width=.65\textwidth]{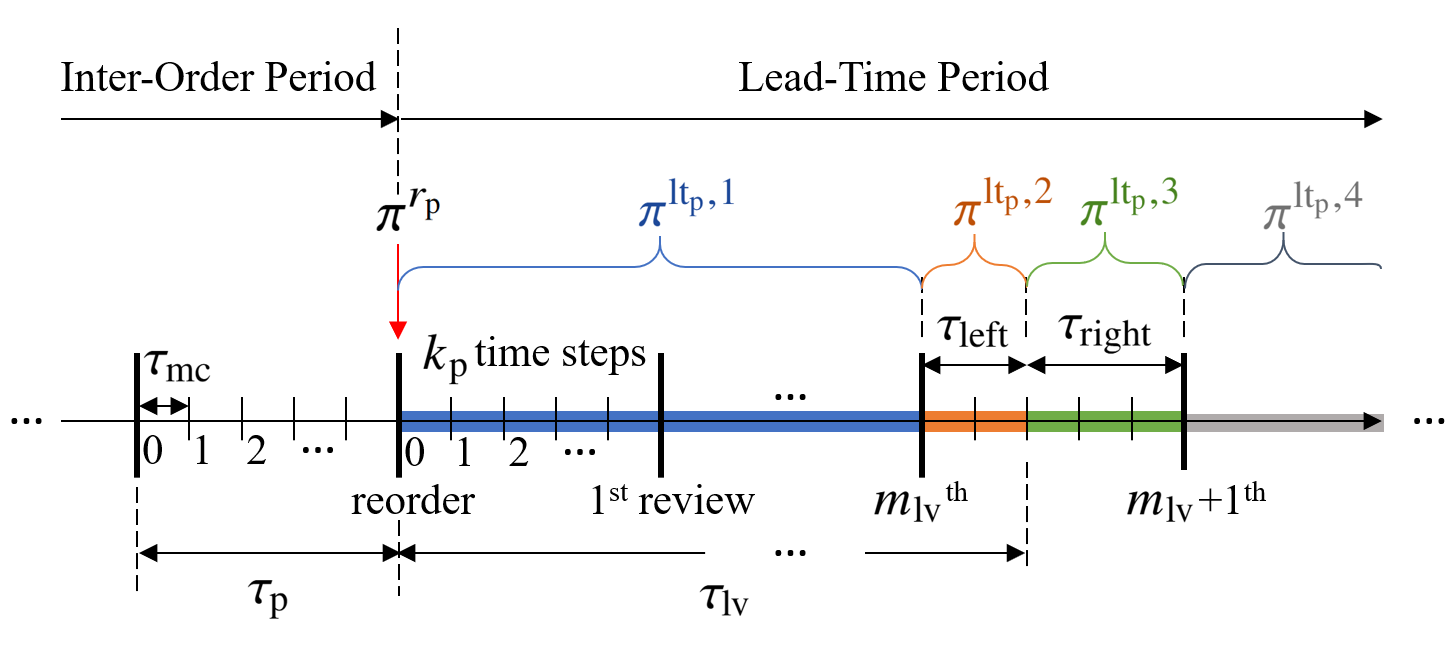}
    \caption{Timeline with the considered lead-time distribution}
    \label{fig:indirect_park_r2q}
\end{figure}

In this work, we can derive a closed-form expression for Eq.\eqref{eq:eta_piq} using the explicit lead-time probability formula from Eq.\eqref{eq:rho_lv}. Referring to Fig.~\ref{fig:indirect_park_r2q}, we define the following integers:
\begin{equation} \label{eq:m_k_def}
    m_\text{lv} = \left\lfloor \frac{\tau_\text{lv}}{\tau_\text{p}} \right\rfloor, \quad 
    k_\text{left} = k_\text{lv} - m_\text{lv} k_\text{p} ,\quad
    k_\text{right} = (m_\text{lv}+1)k_\text{p} - k_\text{lv}
\end{equation}
where $m_\text{lv}$ denotes the index of the first review period during which a delivery can occur, and $k_\text{left}$ indicates the residual number of time steps within the next review period before the fixed processing delay is completed. Here, $\lfloor \cdot \rfloor$ is the floor operator. Since $\tau_{(\cdot)} = k_{(\cdot)} \tau_\text{mc}$, it follows that $\tau_\text{left} = k_\text{left} \tau_\text{mc}$ and $\tau_\text{right} = k_\text{right} \tau_\text{mc}$.

After substituting Eq.~\eqref{eq:rho_lv} into Eq.~\eqref{eq:eta_piq} and careful enumeration, the final expression is
\begin{equation}
    \eta_i \pi^{q_\text{p}}_i = 
    \begin{cases}
        (1 - \alpha) \alpha^{i-1+k_\text{right}} P_{q_\text{p}}\left(P_{f_\text{p}}\right)^{m_\text{lv}+1}\left(I-\alpha^{k_\text{p}} P_{f_\text{p}}  \right)^{-1} \pi^{r_\text{p}}, &\text{if } i \leq k_\text{left} \\
        (1 - \alpha) \alpha^{i-1-k_\text{left}} P_{q_\text{p}}\left(P_{f_\text{p}}\right)^{m_\text{lv}}\left(I-\alpha^{k_\text{p}} P_{f_\text{p}}  \right)^{-1} \pi^{r_\text{p}}, &\text{if } i > k_\text{left}
    \end{cases},
\end{equation}
where $\alpha = e^{-\tau_\text{mc}/\mu_\text{lv}}$, $\forall i=1,\dots,k_\text{p}$.
Similarly, the explicit formula for $\pi^{q_\text{p}}$ in Eq.~\eqref{eq:pi_qp} becomes
\begin{equation} \label{eq:pir_to_piq}
    \pi^{q_\text{p}} = P_{q_\text{p}} \left(P_{f_\text{p}}\right)^{m_\text{lv}} \left( (1 -\alpha^{k_\text{right}}) I + (1-\alpha^{k_\text{p}}) \alpha^{k_\text{right}} P_{f_\text{p}} 
    \left(I - \alpha^{k_\text{p}} P_{f_\text{p}} \right)^{-1} \right)\pi^{r_\text{p}}
\end{equation}

\subsection*{Closed-Form Expression for Eq.~\eqref{eq:pi_lt_p1}}
As before, we can derive a closed form expression for Eq.\eqref{eq:pi_lt_p1}. For the shifted exponential lead-time distribution, the no-delivery probability at step $\ell$ is $\rho_\ell^c = 1$ for $0 \le \ell \le k_\text{lv}$ and $\rho_\ell^c = \alpha^{\ell-k_\text{lv}}$ for $\ell \ge k_\text{lv}+1$. Referring to Fig.\ref{fig:indirect_park_r2q}, the lead time is divided into four segments, each shown in a different color. The distribution $\pi^{\text{lt}_\text{p}}$ is then the weighted sum of the distributions over those four segments. After some algebra and careful enumeration, we obtain
\begin{equation} \label{eq:pi_lt_p2}
    \pi^{\text{lt}_\text{p}} = \frac{1}{k_{\text{lt}_\text{p}}}\sum_{k=1}^{4} \pi^{\text{lt}_\text{p},k}
\end{equation}
where
\begin{equation} \label{eq:pi_lt_p3}
\begin{aligned}
    \pi^{\text{lt}_\text{p},1} &= k_\text{p} \left(I+P_{f_\text{p}} + \cdots + \left(P_{f_\text{p}}\right)^{m_\text{lv}-1} \right)\pi^{r_\text{p}}\\
    \pi^{\text{lt}_\text{p},2} &= (k_\text{left}+1) \left(P_{f_\text{p}}\right)^{m_\text{lv}}\pi^{r_\text{p}}\\
    \pi^{\text{lt}_\text{p},3} &= \frac{ \alpha \cdot ( \alpha^{k_\text{right}-1} - 1) }{ \alpha - 1} \left(P_{f_\text{p}}\right)^{m_\text{lv}} \pi^{r_\text{p}}\\
    \pi^{\text{lt}_\text{p},4} &= \frac{ \alpha^{k_\text{right}} \cdot 
    ( \alpha^{k_\text{p}} - 1) }{ \alpha - 1}
    \left(P_{f_\text{p}}\right)^{m_\text{lv}+1} \left(I - \alpha^{k_\text{p}} P_{f_\text{p}} \right)^{-1}
    \pi^{r_\text{p}}\\
\end{aligned}
\end{equation}
and the unit converted time interval is  $\tau_{\text{lt}_\text{p}} = k_{\text{lt}_\text{p}} \tau_\text{mc}$. 

The first term, $\pi^{\text{lt}_\text{p},1}$, represents the cumulative distribution over the first $m_\text{lv}$ review cycles, each has $k_\text{p}$ time steps. The second term, $\pi^{\text{lt}_\text{p},2}$, corresponds to the $k_\text{left}+1$  residual steps within the $m_\text{lv}+1^\text{th}$ cycle, during which the distribution remains unchanged. The third and fourth terms, $\pi^{\text{lt}_\text{p},3}$ and $\pi^{\text{lt}_\text{p},4}$, are derived from a finite and an infinite geometric series, respectively, under the exponential lead-time model.

\subsection*{Closed-Form Expression for Eq.~\eqref{eq:pi_rc_E}}
Here, we derive the closed-form expressions for $\pi^{\text{io}_\text{p}|E}$ and $\pi^{\text{lt}_\text{p}|E}$ in Eq.~\eqref{eq:pi_rc_E}.  During the IO period, RAAN contact always happens at step $k_\text{p}$ of each review cycle. Average distribution over these points is actually the second term of Eq.~\eqref{eq:pi_io_p_1} as:
\begin{equation}  \label{eq:pi_iop_E}
\begin{aligned}
    \pi^{\text{io}_\text{p}|E} 
    &= \frac{1}{k_{\text{io}_\text{p}|E}} \sum_{j=0}^\infty \left(C^+_{r_\text{p}}P_{f_\text{p}}\right)^j \pi^{q_\text{p}} \\
    &= \frac{1}{k_{\text{io}_\text{p}|E}}\left( I - C^+_{r_\text{p}}P_{f_\text{p}} \right)^{-1} \pi^{q_\text{p}},
\end{aligned}
\end{equation}
where $k_{\text{io}_\text{p}|E}$ is the normalization constant of $\pi^{\text{io}_\text{p}|E}$.

Likewise, during the LT period we apply the same logic to Eq.~\eqref{eq:pi_lt_p1}, again evaluating at step $k_\text{p}$. Enumerating all cases for the LT period gives:
\begin{equation} \label{eq:pi_ltp_E}
\begin{aligned}
    \pi^{\text{lt}_\text{p}|E} 
    &= \frac{1}{k_{\text{lt}_\text{p}|E}}\left[ 
    \sum_{j=1}^{m_\text{lv}} \left( P_{f_\text{p}} \right)^{j-1} + 
    \sum_{j=m_\text{lv}+1}^{\infty} \rho^c_{j k_\text{p}} \left(P_{f_\text{p}}\right)^{j-1} \right]  \pi^{r_\text{p}} \\
    &= \frac{1}{k_{\text{lt}_\text{p}|E}}\left[ 
    \sum_{j=1}^{m_\text{lv}} \left( P_{f_\text{p}} \right)^{j-1}  + \alpha^{k_\text{right}} \left(P_{f_\text{p}}\right)^{m_\text{lv}} \left(I - \alpha^{k_\text{p}} P_{f_\text{p}} \right)^{-1} \right]  \pi^{r_\text{p}}.
\end{aligned}
\end{equation}
where $k_{\text{lt}_\text{p}|E}$ is the normalization constant of $\pi^{\text{lt}_\text{p}|E}$. Referring to Fig.\ref{fig:indirect_park_r2q}, the first term accounts for the cases where $E$ occurs at $j^\text{th}$ review epoch with $j\leq m_\text{lv}$, and the corresponding distribution is $\left(P_{f_\text{p}}\right)^{j-1} \pi^{r_\text{p}}$. The second term represents the cases where $E$ occurs at $j > m_\text{lv}$, with the probability of still awaiting delivery at that time given by $ \rho^c_{j k_\text{p}}$. Note that this expression is sampled only at event $E$, that is, just before parking spares can be distributed to the in-plane orbit, and the result is used to compute the parking availability.

\section*{Appendix B: Direct Strategy as a Special Case of Parking Policy}

The direct strategy is modeled as an $(r,q)$ policy with a shifted exponential lead time distribution~\cite{han2025direct}, while the parking spare policy is modeled as a $(r,q,\tau)$ policy with the same lead time model. If the fixed review period becomes the minimum time step (i.e., $k_\text{p} = 1$), the demand-induced failure matrix becomes the satellite failure matrix (i.e., $P_{f_\text{p}} = P_f$). If both policies share the same $(r,q)$ parameters and the same lead time parameters, they are essentially equivalent by the definition of the $(r,q)$ policy. In this case, we have $k_\text{left} = 0$ and $k_\text{right} = 1$ from Eq.~\eqref{eq:m_k_def} and as shown in Fig.~\ref{fig:indirect_park_r2q}.

Under these conditions, the transition equations between $\pi^q$ and $\pi^r$ in Eqs.~\eqref{eq:piq_to_pir} and \eqref{eq:pi_qp} become:
\begin{equation}
\begin{aligned}
    \pi^r &= C_r^- P_f (I - C_r^+ P_f)^{-1} \pi^q, \\
    \pi^q &= \sum_{j=0}^\infty \rho_{j+1} P_q \left(P_f\right)^j \pi^r \\
          &= (1 - \alpha) P_q \left(P_f\right)^{m} \left(I - \alpha P_f\right)^{-1} \pi^r.
\end{aligned}
\end{equation}
Likewise, the distributions during the IO and LT periods  in Eqs.~\eqref{eq:pi_io_p_2} and \eqref{eq:pi_lt_p1}  are:
\begin{equation}
\begin{aligned}
    \pi^{\text{io}} &= \frac{1}{k_{\text{io}}} C_r^+ P_f (I - C_r^+ P_f)^{-1} \pi^q, \\
    \pi^{\text{lt}} &= \frac{1}{k_{\text{lt}_\text{p}}} \sum_{k=1}^{4} \pi^{\text{lt}_\text{p},k} \\
    &= \frac{1}{k_{\text{lt}}} \left( \sum_{i=0}^m \left(P_f\right)^i + \alpha \left(P_f\right)^{m+1} \left(I - \alpha P_f\right)^{-1} \right) \pi^r.
\end{aligned}
\end{equation}
All four expressions exactly match the results presented in Ref.~\cite{han2025direct}. Thus, this paper generalizes the previous work and develops a unified formulation that covers both the direct and indirect strategies.

Using the above equations, one can compute $\pi^{\text{rc}}$ and $\tau_{\text{rc}}$ as described in Sec.~\ref{sec:park_rc}. With all stationary distributions computed, the cost and resilience metrics of the direct strategy can be obtained by omitting the parking-orbit related terms in Sec.~\ref{sec4} (e.g., parking spare holding cost and orbit transfer cost).

\section*{Funding Sources}
This research was supported by the Advanced Technology R\&D Center at Mitsubishi Electric Corporation.

\section*{Acknowledgments}
Artificial intelligence tools were used to improve the grammar and readability of this paper and to assist with reference formatting. The authors reviewed and verified the final content.

\bibliography{references}

\end{document}